\documentclass[10pt,reqno]{amsart}%
\usepackage[T1]{fontenc}
\usepackage{xcolor}
\PassOptionsToPackage{hyphens}{url}
\usepackage{hyperref}
\usepackage{amsfonts}
\usepackage{amssymb}
\usepackage{amsthm}
\usepackage{mathtools}
\usepackage{mathptmx}
\usepackage{geometry}
\usepackage{enumitem}
\setlist[enumerate]{leftmargin=1.6em,labelsep=0.4em,itemsep=2pt,topsep=4pt}
\setlist[itemize]{leftmargin=1.6em,labelsep=0.4em,itemsep=2pt,topsep=4pt}
\usepackage[numbers,sort&compress]{natbib}
\usepackage[expansion=false]{microtype}
\usepackage{graphicx}%
\usepackage{amsmath}%
\hypersetup{
pdftitle={Associative operators and their iterated lifts: classification and cocycle identities},
pdfauthor={Dorian Popa, Mircea Dan Rus and Constantin-Cosmin Todea},
colorlinks=true,
linkcolor={blue!60!black},
citecolor={blue!60!black},
urlcolor={blue!80!black},
}
\providecommand{\doi}{}
\renewcommand{\doi}[1]{\href{https://doi.org/#1}{\nolinkurl{https://doi.org/#1}}}
\theoremstyle{plain}
\newtheorem{theorem}{Theorem}[section]
\newtheorem{proposition}[theorem]{Proposition}
\newtheorem{lemma}[theorem]{Lemma}
\newtheorem{corollary}[theorem]{Corollary}
\theoremstyle{definition}
\newtheorem{definition}[theorem]{Definition}
\newtheorem{example}[theorem]{Example}
\newtheorem{remark}[theorem]{Remark}
\newtheorem{question}[theorem]{Question}
\begin{document}
\title[Associative operators and their iterated lifts]{Associative operators and their iterated lifts:\\ classification and cocycle identities}

\author[D. Popa]{Dorian Popa}
\address{Department of Mathematics, Technical University of Cluj-Napoca, 28
Memorandumului Street, 400114 Cluj-Napoca, Romania\\
ORCID: \href{https://orcid.org/0000-0001-6197-3138}{0000-0001-6197-3138}}
\email{popa.dorian@math.utcluj.ro}

\author[M.~D. Rus]{Mircea Dan Rus}
\address{Department of Mathematics, Technical University of Cluj-Napoca, 28
Memorandumului Street, 400114 Cluj-Napoca, Romania\\
ORCID: \href{https://orcid.org/0000-0001-6255-0241}{0000-0001-6255-0241}}
\email{rus.mircea@math.utcluj.ro}

\author[C.-C. Todea]{Constantin-Cosmin Todea}
\address{Department of Mathematics, Technical University of Cluj-Napoca, 28
Memorandumului Street, 400114 Cluj-Napoca, Romania\\
ORCID: \href{https://orcid.org/0000-0002-0277-1285}{0000-0002-0277-1285}}
\email{Constantin.Todea@math.utcluj.ro}

\date{11 November 2026}
\subjclass[2020]{Primary 39B52; Secondary 39A70, 20J06, 20M30, 20M50}
\keywords{functional equation, Cauchy difference, quadratic functional equation, associative operator, tower, cocycle identity, inverse problem, symmetry}

\begin{abstract}
The inverse problem of an inhomogeneous functional equation asks for which right-hand sides the equation is solvable. For the Cauchy equation and for the quadratic equation a necessary condition is a cocycle identity, and the two identities have a common form. We take this form as a definition. An operator from functions of one variable to functions of two variables can be lifted to functions of several variables, one argument at a time. We call it associative when its two iterated lifts from one variable to three coincide, that is, when every function in its image satisfies the corresponding cocycle identity. Composition with an associative binary operation is the model case. The Cauchy difference, the d'Alembert operator and the Jordan--von Neumann operator are associative as well. We prove that all iterated lifts of the same length of an associative operator coincide. Over abelian groups we classify, under suitable hypotheses, the associative integer combinations of translations, embeddings and the evaluation at the origin. For a composition operator the cocycle identity characterises the image over monoids and over bands. Over an arbitrary semigroup we reduce this question to the injectivity of a canonical map.
\end{abstract}

\maketitle

\section{Introduction}\label{sec:intro}

Throughout the paper $G$ and $H$ are non-empty sets. For $n\ge 1$ we write $\mathcal{H}_n := H^{G^{n}}$ for the set of all functions $G^n\to H$, and we call the sequence $(\mathcal{H}_n)_{n\ge 1}$ the \emph{function hierarchy} of the pair $(G,H)$. An operator $\mathsf{T}\colon \mathcal{H}_1\to \mathcal{H}_2$ turns functions of one variable into functions of two. It then acts on the whole hierarchy, one argument at a time, through the lifts $\mathsf{T}^{[i]}\colon \mathcal{H}_n\to \mathcal{H}_{n+1}$ of Section~\ref{sec:lifts}, where the index $i$ is the position, or \emph{slot}, of the argument on which $\mathsf{T}$ acts. There are two ways of iterating $\mathsf{T}$ from $\mathcal{H}_1$ to $\mathcal{H}_3$, namely $\mathsf{T}^{[1]}\mathsf{T}$ and $\mathsf{T}^{[2]}\mathsf{T}$, and we call $\mathsf{T}$ \emph{associative} when these two composites coincide. The model case is composition with a binary operation $\mu$ on $G$, that is, the operator $\mathsf{T}[f](x,y) = f\bigl(\mu(x,y)\bigr)$. This operator is associative whenever $\mu$ is, and the converse holds as soon as $H$ has at least two elements (Proposition~\ref{prop:model}).

The motivation comes from the \emph{inverse problems} of inhomogeneous functional equations: given an operator $\mathsf{T}$ and a datum $\psi\in \mathcal{H}_2$, decide when the equation $\mathsf{T}[f]=\psi$ is solvable. This question is also called the \emph{image problem} for~$\mathsf{T}$. When the equation is solvable, one also wants to determine $f$ from~$\psi$. For the Cauchy difference $\mathsf{L}[f](x,y) = f(x+y)-f(x)-f(y)$, whose homogeneous equation is the subject of the monograph \cite{Kuczma}, the classical answer describes the image as the set of symmetric solutions of the cocycle equation, under hypotheses that Section~\ref{sec:image} states. That line of results goes from Kurepa and Erdős \cite{Kurepa1956,Erdos1959} to Davison and Ebanks \cite{DavisonEbanks1995}, and Section~\ref{sec:image} recalls it with its literature. The inhomogeneous equation is studied in \cite{FenyoForti,AlievIsmailov}. For the Jordan--von Neumann operator $\mathsf{Q}[f](x,y)=f(x+y)+f(x-y)-2f(x)-2f(y)$ of the quadratic functional equation, a cocycle identity is again a necessary condition on the datum, and the inverse problem is studied in \cite{PopaRusJvN}. The two identities are derived separately, and \cite[Section~2]{PopaRusJvN} observes that they have a common two-slot form. We take that form as a definition. An operator $\mathsf{T}$ is associative exactly when every function $\psi$ in its image satisfies $\mathsf{T}^{[1]}[\psi]=\mathsf{T}^{[2]}[\psi]$ (Proposition~\ref{prop:cocycle-complex}), and for $\mathsf{L}$ and for $\mathsf{Q}$ this is the identity known from the literature. The composition operators of associative operations are associative as well. A similar identity, in four variables, characterises the solvable data of the inhomogeneous equation of Sincov \cite{MadutaPopaSingeorzan}.

The inverse problems are also of interest for the stability theory of functional equations, which goes back to Hyers's answer \cite{Hyers1941} to a question of Ulam and is presented in the monograph \cite{HyersIsacRassias}. When an inverse theorem gives the solution by an explicit formula in terms of the datum, a bound on the datum can yield a bound on the solution.

Composition with an associative operation is the model case, but not the only one. Over abelian groups, the Cauchy difference is associative, and so are the d'Alembert operator $\mathsf{P}[f](x,y)=f(x+y)+f(x-y)$ and the Jordan--von Neumann operator~$\mathsf{Q}$, while the plain finite difference $\Delta[f](x,y) = f(x+y)-f(y)$ is not associative unless $G$ or $H$ is trivial (Proposition~\ref{prop:delta-not}). We call an iterated lift $\mathsf{T}^{[j_n]}\cdots\mathsf{T}^{[j_2]}\mathsf{T}^{[j_1]}\colon \mathcal{H}_1\to \mathcal{H}_{n+1}$ a \emph{tower} of length~$n$ over~$\mathsf{T}$. Our first result, the tower theorem (Theorem~\ref{thm:towers}), asserts that all towers of the same length over an associative operator coincide, for arbitrary $G$ and $H$. Their common value is the power $\mathsf{T}^{\,n}\colon \mathcal{H}_1\to \mathcal{H}_{n+1}$. Associativity does not imply the exchange law $\mathsf{T}^{[i]}\mathsf{T}^{[j]}=\mathsf{T}^{[j+1]}\mathsf{T}^{[i]}$, $i<j$, between lifts in distinct slots (Example~\ref{ex:no-exchange}), so the lifts of an associative operator need not satisfy the coface relations of a cosimplicial object. For $\mathsf{L}$ the powers are the higher-order Cauchy differences. For the non-associative $\Delta$ we describe all the towers. They are mixed differences in some of the variables, evaluated at the sum of the others, and we call them partial mixed differences (Section~\ref{sec:towers}). The same section computes the powers of the Jordan--von Neumann operator in closed form, and identifies the mixed differences at the origin, from the second one on, as the powers of a single associative operator. The numbers of distinct towers of an operator form its tower spectrum, which is bounded by $n!$ at length~$n$, and the bound is attained (Proposition~\ref{prop:tower-spectrum}).

The remaining results are as follows. Section~\ref{sec:sums} determines, in a single classification over abelian groups (Theorem~\ref{thm:global}), the associative operators in the five-parameter integer family generated by the two translations, the two embeddings and the evaluation at the origin. Its hypothesis on the target cannot be dropped, and its hypothesis on the domain serves only the criterion in terms of the coefficients (Remark~\ref{rem:sharp}). The classification characterises, up to integer multiples, the Cauchy difference, the d'Alembert and Jordan--von Neumann operators and the mixed difference at the origin inside this family (Corollary~\ref{cor:detect}). Section~\ref{sec:symmetry} studies symmetry. For a target with at least two elements, a composition operator is symmetric exactly when its operation is commutative (Proposition~\ref{prop:symmetric-model}). The powers of a symmetric associative operator take symmetric functions as values (Theorem~\ref{thm:symmetric-towers}), and within the classification the symmetric associative operators are exactly the members of the two pencils $\{m\mathsf{L}+n\mathsf{J}_0\}$ and $\{m\mathsf{E}_{+}+n(\mathsf{J}_0-\mathsf{J}_1-\mathsf{J}_2)\}$ of Theorem~\ref{thm:global} (Proposition~\ref{prop:symmetric-classification}). Section~\ref{sec:image} returns to the inverse problems. It states the cocycle identities with their instances and shows that for a symmetric datum they reduce to a single identity. For a composition operator it characterises the image over a monoid and over a band, and it reduces the same question over an arbitrary semigroup to the injectivity of a canonical map (Proposition~\ref{prop:balanced}). The weighted composition operators over monoids and over bands are reduced to the unweighted ones. Section~\ref{sec:new} separates what is new from what is known, with the corresponding attributions, and states three open questions.

\subsection*{Notation and conventions}
Operators are set in sans-serif and applied with brackets, $\mathsf{T}[f]$. The composition operators $\mu^{*}$ and their weighted forms $\omega\mu^{*}$ keep the letter of their operation. The finite differences $\Delta$ and $\Delta_t$ keep the upright capital delta traditional in the difference calculus. Juxtaposition of operators denotes composition, and a product $\prod$ of commuting operators denotes their composite. A superscript in square brackets is a slot, as in the lift $\mathsf{T}^{[i]}$, and a plain superscript is a power, as in $\mathsf{T}^{\,n}$. For an operator $\mathcal{H}_1\to\mathcal{H}_1$, such as $\Delta_t$, a plain superscript is the usual power under composition. Subscripts follow no single convention, and each is fixed where the operator carrying it is introduced. Those of the embeddings $\mathsf{J}_1$ and $\mathsf{J}_2$ are slot indices, those of the evaluations $\mathsf{J}_e$ and of the shifts $\mathsf{E}_t$ and differences $\Delta_t$ are elements of~$G$, and the signs of $\mathsf{E}_{+}$ and $\mathsf{E}_{-}$ are neither. The symbols $\mathsf{J}_1$ and $\mathsf{J}_2$ always denote the embeddings, whatever the elements of $G$ happen to be called.

The letters $f$ and $u$ are kept for elements of $\mathcal{H}_1$, $\psi$ and $\chi$ for elements of a general $\mathcal{H}_n$, $x,y,z$ and the indexed $x_1,x_2,\dots$ for elements of $G$, and $c$ for elements of~$H$. When a constant $c$ is equated with the value of an operator, it stands for the corresponding constant function. A centred dot in an argument list marks the place of the variable of a partial map, as in $\psi(x,\,\cdot\,,z)$. The letters $i$ and $j$ index slots, $k$ and $\ell$ index summands and elements of index sets, and $n$ is an arity or the length of a tower. The letters $m$ and $n$ also serve as integer coefficients, as in the pencil $\{m\mathsf{L}+n\mathsf{J}_0\}_{m,n\in\mathbb{Z}}$, where the subscript records the use. All other letters are introduced where they occur, and keep their meaning only in the statement or proof that introduces them. We write $\mathbb{Z}$ for the ring of integers and $\mathbb{Z}_m$ for the integers modulo~$m$. For an integer $m$, the notation $mG=0$ means that $mg=0$ for every $g\in G$, and similarly for~$H$. All rings are commutative with identity, and all modules are unital.

Any structure on $G$ or on $H$ beyond that of a non-empty set is stated with the result that uses it. A semigroup operation on $G$ is written additively, and the neutral element of a monoid is~$0$. Additive notation usually signals a commutative operation. We keep it for arbitrary semigroups, because the Cauchy difference and the related operators are traditionally written with $x+y$, and we assume commutativity only where this is stated. When the operation is a general map $\mu\colon G^2\to G$, as for the composition operators, it is written $\mu(x,y)$, and the identity element of a monoid is then denoted by~$e$. The operation of $H$, where $H$ carries one, is written additively as well, and the two operations are distinguished by their arguments. A general map $H^2\to H$ is written~$\beta$, as in Example~\ref{ex:no-exchange}. When $H$ is a semigroup, each $\mathcal{H}_n$ is a semigroup under the pointwise operation, sums of operators are formed pointwise, and an operator between such sets is called \emph{additive} when it preserves that operation, that is, when it is a homomorphism for the pointwise operations. Commutativity of $H$ is assumed where the order of the terms matters, as in Proposition~\ref{prop:P}. When $H$ is an abelian group, so is each $\mathcal{H}_n$, and integer combinations of operators are formed pointwise.

\section{Lifts, associativity and first examples}\label{sec:lifts}

\begin{definition}\label{def:lift}
Let $\mathsf{T}\colon \mathcal{H}_1\to \mathcal{H}_2$ be an operator, and let $n\ge 1$ and $1\le i\le n$. The operator $\mathsf{T}^{[i]}\colon \mathcal{H}_n\to \mathcal{H}_{n+1}$ defined by
\begin{equation}\label{eq:lift}
\mathsf{T}^{[i]}[\psi](x_1,\dots,x_{n+1}) := \mathsf{T}\bigl[\psi(x_1,\dots,x_{i-1},\,\cdot\,,x_{i+2},\dots,x_{n+1})\bigr](x_i,x_{i+1})
\end{equation}
is called the \emph{lift of $\mathsf{T}$ in the $i$-th slot}.
\end{definition}

The $i$-th variable of $\psi$ is replaced by the pair $(x_i,x_{i+1})$, because $\mathsf{T}$ produces functions of two variables. We call $(x_i,x_{i+1})$ the \emph{active pair} of the lift, and the remaining variables, which keep their relative order, \emph{passive}. On $\mathcal{H}_1$ there is a single slot and $\mathsf{T}^{[1]}=\mathsf{T}$ there. The same symbol $\mathsf{T}^{[i]}$ denotes the lift at every arity $n\ge i$, and the argument determines the arity. This is the usage familiar from the coface maps of simplicial theory, and it leaves the subscript position free for operators that carry indices of their own, such as the embeddings $\mathsf{J}_1,\mathsf{J}_2$ of Example~\ref{ex:comp}.

\begin{definition}\label{def:composition}
Let $\mu\colon G^2\to G$ be a map. The operator $\mu^{*}\colon \mathcal{H}_1\to \mathcal{H}_2$, $\mu^{*}[f] := f\circ \mu$, is called the \emph{composition operator} of~$\mu$. More generally, for $r,s\ge 1$ and a map $\nu\colon G^{s}\to G^{r}$ we write $\nu^{*}\colon \mathcal{H}_r\to \mathcal{H}_s$ for the operator $\nu^{*}[\psi] := \psi\circ\nu$. For $n\ge 1$ and $1\le i\le n$, the map $\mu_{[i]}\colon G^{n+1}\to G^{n}$ that applies $\mu$ to the slots $i$ and $i+1$ is defined by
\[
\mu_{[i]}(x_1,\dots,x_{n+1}) := \bigl(x_1,\dots,x_{i-1},\mu(x_i,x_{i+1}),x_{i+2},\dots,x_{n+1}\bigr).
\]
\end{definition}

For every $\psi\in\mathcal{H}_n$, composition with $\mu_{[i]}$ substitutes $\mu(x_i,x_{i+1})$ for the $i$-th argument of $\psi$ and keeps the order of the others,
\[
(\psi\circ\mu_{[i]})(x_1,\dots,x_{n+1}) = \psi\bigl(x_1,\dots,x_{i-1},\mu(x_i,x_{i+1}),x_{i+2},\dots,x_{n+1}\bigr).
\]

We record the elementary properties of the lifts.

\begin{lemma}\label{lem:basic}
Let $n\ge 1$ and $1\le i\le n$.
\begin{enumerate}
\item[\rm (i)] The lift of a composition operator is a composition operator: $(\mu^{*})^{[i]} = (\mu_{[i]})^{*}$ for every map $\mu\colon G^2\to G$.
\item[\rm (ii)] If $H$ is a semigroup and $\mathsf{T},\mathsf{S}\colon \mathcal{H}_1\to \mathcal{H}_2$, then $(\mathsf{T}+\mathsf{S})^{[i]} = \mathsf{T}^{[i]} + \mathsf{S}^{[i]}$.
\item[\rm (iii)] If $H$ is a semigroup and $\mathsf{T}\colon \mathcal{H}_1\to \mathcal{H}_2$ is additive, then $\mathsf{T}^{[i]}$ is additive.
\end{enumerate}
\end{lemma}

\begin{proof}
All three statements are immediate from~(\ref{eq:lift}). For~(i), applying $\mu^{*}$ to the $i$-th partial map of $\psi$ and evaluating at $(x_i,x_{i+1})$ returns the value of $\psi$ at the tuple $\mu_{[i]}(x_1,\dots,x_{n+1})$, so $(\mu^{*})^{[i]}[\psi] = \psi\circ \mu_{[i]}$. The lift is additive in the operator in~(ii) because $(\mathsf{T}+\mathsf{S})[u]=\mathsf{T}[u]+\mathsf{S}[u]$ for every $u\in \mathcal{H}_1$, and it is additive in the argument in~(iii) because passing to a partial map is additive.
\end{proof}

\begin{definition}\label{def:associative}
An operator $\mathsf{T}\colon \mathcal{H}_1\to \mathcal{H}_2$ is said to be \emph{associative} if $\mathsf{T}^{[1]}\mathsf{T} = \mathsf{T}^{[2]}\mathsf{T}$. The common value is denoted $\mathsf{T}^{\,2}\colon \mathcal{H}_1\to \mathcal{H}_3$.
\end{definition}

Definition~\ref{def:associative} has a precursor in \cite[Section~2]{PopaRusJvN}, where the cocycle identities of the Cauchy difference and of the Jordan--von Neumann operator are observed to share a common two-slot form. New here are the notion itself and the lifts in every slot, through which the operator acts on the whole hierarchy.

The model case explains the name. When an operation $\mu$ on $G$ is associative, we write $\mu(x,y,z)$ for the common value of $\mu(\mu(x,y),z)$ and $\mu(x,\mu(y,z))$. More generally, for $n\ge 2$ we write $\mu(x_1,\dots,x_n)$ for the common value of all the bracketings of $x_1,\dots,x_n$, which coincide by the generalised associative law \cite[Ch.~I, §1, no.~3]{Bourbaki}, so that
\[
\mu(x_1,\dots,x_{n+1}) = \mu\bigl(\mu(x_1,\dots,x_n),x_{n+1}\bigr) = \mu\bigl(x_1,\mu(x_2,\dots,x_{n+1})\bigr),
\]
with the convention $\mu(x)=x$.

\begin{proposition}\label{prop:model}
Let $\mu\colon G^2\to G$ be a map. If $\mu$ is associative, then the composition operator $\mu^{*}$ is associative and $(\mu^{*})^{\,2}[f](x,y,z) = f\bigl(\mu(x,y,z)\bigr)$. If $H$ has at least two elements, the converse holds as well.
\end{proposition}

\begin{proof}
By Lemma~\ref{lem:basic}(i), $(\mu^{*})^{[1]}\mu^{*}[f] = (\mu_{[1]})^{*}\bigl[f\circ\mu\bigr] = f\circ\mu\circ\mu_{[1]}$ and, in the same way, $(\mu^{*})^{[2]}\mu^{*}[f] = f\circ\mu\circ\mu_{[2]}$, that is,
\[
(\mu^{*})^{[1]}\mu^{*}[f](x,y,z) = f\bigl(\mu(\mu(x,y),z)\bigr), \qquad (\mu^{*})^{[2]}\mu^{*}[f](x,y,z) = f\bigl(\mu(x,\mu(y,z))\bigr).
\]
If $\mu$ is associative the two agree for every~$f$, with common value $f(\mu(x,y,z))$. Conversely, if $\mu(\mu(x,y),z)\neq \mu(x,\mu(y,z))$ for some $x,y,z$, choose $c_0\neq c_1$ in~$H$ and let $f$ take the value $c_0$ at $\mu(\mu(x,y),z)$ and $c_1$ elsewhere. The two sides then differ at $(x,y,z)$.
\end{proof}

Example~\ref{ex:comp} lists the associative operations used in the paper, each with the composition operator it yields, which is associative by Proposition~\ref{prop:model}.

\begin{example}[composition operators]\label{ex:comp}
\leavevmode
\begin{enumerate}
\item[\rm (i)] \emph{The translation.} Let $G$ be a semigroup. The addition, being associative, yields the associative operator $\mathsf{E}[f](x,y) := f(x+y)$, with $\mathsf{E}^{\,2}[f](x,y,z) = f(x+y+z)$.
\item[\rm (ii)] \emph{The embeddings.} Over an arbitrary set $G$, the two projections are associative operations and yield the associative operators $\mathsf{J}_1[f](x,y) := f(x)$ and $\mathsf{J}_2[f](x,y) := f(y)$, with $\mathsf{J}_1^{\,2}[f](x,y,z) = f(x)$ and $\mathsf{J}_2^{\,2}[f](x,y,z) = f(z)$.
\item[\rm (iii)] \emph{Evaluations.} For a fixed $e\in G$, the constant operation $\mu(x,y)=e$ is associative and yields the operator $\mathsf{J}_e[f](x,y) := f(e)$, with $\mathsf{J}_e^{\,2}[f]=f(e)$. When $G$ is a monoid we write $\mathsf{J}_0$ for the operator obtained from $e=0$, that is, $\mathsf{J}_0[f](x,y) := f(0)$, and call it the \emph{evaluation at the origin}.
\item[\rm (iv)] \emph{Semilattices and bands.} Over a linearly ordered set~$G$, the operation $\max$ is associative, and the composition operator $(\max^{*})[f](x,y)=f(\max(x,y))$ is an associative operator, with $(\max^{*})^{\,2}[f](x,y,z)=f(\max(x,y,z))$; similarly for $\min$. The operation $\max$ is a special case of the following. A semigroup in which every element is idempotent is called a \emph{band} \cite[§4.4]{Howie}, and a commutative band is called a \emph{semilattice}. The operation of a band is associative by definition, so its composition operator is an associative operator, again by Proposition~\ref{prop:model}. A linearly ordered set under $\max$ is a semilattice, and so is any lattice under its join. Bands return in Proposition~\ref{prop:bands}.
\end{enumerate}
\end{example}

The main operators of this paper are integer combinations of composition operators. Whether such a combination is associative cannot be read off from an operation. The first two are the \emph{finite difference} and the \emph{Cauchy difference}: for $G$ a semigroup and $H$ an abelian group, set
\begin{equation*}
\Delta := \mathsf{E}-\mathsf{J}_2, \qquad \mathsf{L} := \mathsf{E}-\mathsf{J}_1-\mathsf{J}_2 = \Delta - \mathsf{J}_1,
\end{equation*}
that is, $\Delta[f](x,y) = f(x+y)-f(y)$ and $\mathsf{L}[f](x,y) = f(x+y)-f(x)-f(y)$. In both definitions the operation enters only through $x+y$, so a semigroup suffices.

\begin{proposition}\label{prop:delta-not}
Let $G$ be a semigroup with $x_0+z_0\neq z_0$ for some $x_0,z_0\in G$, and let $H$ be a non-trivial abelian group. Then the finite difference $\Delta$ is not associative. Precisely, for all $f\in \mathcal{H}_1$ and $x,y,z\in G$,
\begin{equation}\label{eq:delta-failure}
\Delta^{[1]}\Delta[f](x,y,z) - \Delta^{[2]}\Delta[f](x,y,z) = \Delta[f](x,z).
\end{equation}
\end{proposition}

\begin{proof}
From (\ref{eq:lift}), $\Delta^{[1]}[\psi](x,y,z) = \psi(x+y,z)-\psi(y,z)$ and $\Delta^{[2]}[\psi](x,y,z) = \psi(x,y+z)-\psi(x,z)$. With $\psi=\Delta[f]$,
\begin{align*}
\Delta^{[1]}\Delta[f](x,y,z) &= \bigl[f(x+y+z)-f(z)\bigr]-\bigl[f(y+z)-f(z)\bigr] = f(x+y+z)-f(y+z),\\
\Delta^{[2]}\Delta[f](x,y,z) &= \bigl[f(x+y+z)-f(y+z)\bigr]-\bigl[f(x+z)-f(z)\bigr],
\end{align*}
and the difference of the two lines is $f(x+z)-f(z) = \Delta[f](x,z)$. For the function $f$ with value $c\neq 0$ at $x_0+z_0$ and $0$ elsewhere, this is non-zero at $(x_0,y,z_0)$.
\end{proof}

\begin{proposition}\label{prop:L}
Let $G$ be a semigroup and $H$ an abelian group. The Cauchy difference is associative, with $\mathsf{L}^{\,2}=\mathsf{F}$, the Fréchet operator
\[
\mathsf{F}[f](x,y,z) = f(x+y+z)-f(x+y)-f(y+z)-f(x+z)+f(x)+f(y)+f(z).
\]
\end{proposition}

\begin{proof}
With $\psi=\mathsf{L}[f]$, using $\mathsf{L}^{[1]}[\psi](x,y,z)=\psi(x+y,z)-\psi(x,z)-\psi(y,z)$,
\begin{align*}
\mathsf{L}^{[1]}\mathsf{L}[f](x,y,z) &= \bigl[f(x+y+z)-f(x+y)-f(z)\bigr] - \bigl[f(x+z)-f(x)-f(z)\bigr]\\
&\qquad - \bigl[f(y+z)-f(y)-f(z)\bigr] = \mathsf{F}[f](x,y,z),
\end{align*}
and likewise, using $\mathsf{L}^{[2]}[\psi](x,y,z)=\psi(x,y+z)-\psi(x,y)-\psi(x,z)$,
\begin{align*}
\mathsf{L}^{[2]}\mathsf{L}[f](x,y,z) &= \bigl[f(x+y+z)-f(x)-f(y+z)\bigr] - \bigl[f(x+y)-f(x)-f(y)\bigr]\\
&\qquad - \bigl[f(x+z)-f(x)-f(z)\bigr] = \mathsf{F}[f](x,y,z). \qedhere
\end{align*}
\end{proof}

Propositions \ref{prop:delta-not} and~\ref{prop:L} are our starting point. The difference $\Delta$ and the Cauchy difference $\mathsf{L}$ are both obtained from the translation~$\mathsf{E}$ by subtracting embeddings, yet $\mathsf{L}$ is associative and, under the hypotheses of Proposition~\ref{prop:delta-not}, $\Delta$ is not.

The next two operators are of a different kind. Unlike $\Delta$ and $\mathsf{L}$, they involve the inversion of the domain, through the second translation
\[
\mathsf{E}_{-}[f](x,y) := f(x-y),
\]
defined over an abelian group $G$. We write $\mathsf{E}_{+}$ instead of $\mathsf{E}$ when the two are used together. Subtraction is associative only when $2G=0$, since $(x-y)-z$ and $x-(y-z)$ differ by $2z$. By Proposition~\ref{prop:model}, applied to a target $H$ with at least two elements, the operator $\mathsf{E}_{-}$ is therefore not associative unless $2G=0$. The two operators built from it are nevertheless associative. The first is $\mathsf{P}=\mathsf{E}_{+}+\mathsf{E}_{-}$,
\[
\mathsf{P}[f](x,y) := f(x+y)+f(x-y).
\]
We call it the \emph{d'Alembert operator}, after d'Alembert's functional equation $f(x+y)+f(x-y)=2f(x)f(y)$ \cite{AczelDhombres,Stetkaer}, whose left-hand side it is. The equation is classically studied for complex-valued functions on groups, abelian or not \cite{Stetkaer}. The operator $\mathsf{P}$ is considered here on an abelian group, with values in a commutative semigroup.

\begin{proposition}\label{prop:P}
Let $G$ be an abelian group and $H$ a commutative semigroup. The operator $\mathsf{P}$ is associative, with
\[
\mathsf{P}^{\,2}[f](x,y,z) = f(x+y+z)+f(x+y-z)+f(x-y+z)+f(x-y-z).
\]
\end{proposition}

\begin{proof}
With $\psi=\mathsf{P}[f]$,
\begin{align*}
\mathsf{P}^{[1]}\mathsf{P}[f](x,y,z) &= \psi(x+y,z)+\psi(x-y,z)\\
&= f(x+y+z)+f(x+y-z)+f(x-y+z)+f(x-y-z),
\end{align*}
and $\mathsf{P}^{[2]}\mathsf{P}[f](x,y,z) = \psi(x,y+z)+\psi(x,y-z)$ expands to the same four terms, in the order $f(x+y+z)$, $f(x-y-z)$, $f(x+y-z)$, $f(x-y+z)$.
\end{proof}

The second is the Jordan--von Neumann operator $\mathsf{Q}=\mathsf{P}-2\mathsf{J}_1-2\mathsf{J}_2$,
\[
\mathsf{Q}[f](x,y) := f(x+y)+f(x-y)-2f(x)-2f(y),
\]
defined for abelian groups $G$ and $H$. The quadratic functional equation $f(x+y)+f(x-y)=2f(x)+2f(y)$ is the equation $\mathsf{Q}[f]=0$. It goes back to the characterisation of inner product spaces by Jordan and von Neumann \cite{JordanVonNeumann1935}. When the target is uniquely $2$-divisible, that is, when multiplication by $2$ is a bijection of the target, the solutions of the quadratic equation are the diagonals $x\mapsto B(x,x)$ of the symmetric biadditive maps $B$ \cite{Aczel1965}. Up to sign, $\mathsf{Q}$ is the \emph{quadratic difference operator} of \cite{CzerwikDlutek}, defined there as $2f(x)+2f(y)-f(x+y)-f(x-y)$. By Theorem~\ref{thm:global}, under its hypotheses on $G$ and $H$, the integer multiples of $\mathsf{P}$ and of $\mathsf{Q}$ are the only associative members of the family considered there in which $\mathsf{E}_{-}$ has a non-zero coefficient.

\begin{proposition}\label{prop:Q}
Let $G$ and $H$ be abelian groups. The operator $\mathsf{Q}$ is associative, with
\begin{multline}\label{eq:Q2}
\mathsf{Q}^{\,2}[f](x,y,z) = f(x+y+z)+f(x+y-z)+f(x-y+z)+f(x-y-z)\\
-2\bigl[f(x+y)+f(x-y)+f(x+z)+f(x-z)+f(y+z)+f(y-z)\bigr]\\
+4\bigl[f(x)+f(y)+f(z)\bigr].
\end{multline}
\end{proposition}

\begin{proof}
With $\psi=\mathsf{Q}[f]$,
\begin{align*}
\mathsf{Q}^{[1]}\mathsf{Q}[f](x,y,z) &= \psi(x+y,z)+\psi(x-y,z)-2\psi(x,z)-2\psi(y,z)\\
&= \bigl[f(x+y+z)+f(x+y-z)-2f(x+y)-2f(z)\bigr]\\
&\quad + \bigl[f(x-y+z)+f(x-y-z)-2f(x-y)-2f(z)\bigr]\\
&\quad - 2\bigl[f(x+z)+f(x-z)-2f(x)-2f(z)\bigr]\\
&\quad - 2\bigl[f(y+z)+f(y-z)-2f(y)-2f(z)\bigr],
\end{align*}
which collects to the right-hand side of~(\ref{eq:Q2}). Likewise
\begin{align*}
\mathsf{Q}^{[2]}\mathsf{Q}[f](x,y,z) &= \psi(x,y+z)+\psi(x,y-z)-2\psi(x,y)-2\psi(x,z)\\
&= \bigl[f(x+y+z)+f(x-y-z)-2f(x)-2f(y+z)\bigr]\\
&\quad + \bigl[f(x+y-z)+f(x-y+z)-2f(x)-2f(y-z)\bigr]\\
&\quad - 2\bigl[f(x+y)+f(x-y)-2f(x)-2f(y)\bigr]\\
&\quad - 2\bigl[f(x+z)+f(x-z)-2f(x)-2f(z)\bigr],
\end{align*}
which collects to the same expression.
\end{proof}

In the language of Definition~\ref{def:associative}, the associativity of $\mathsf{Q}$ is \cite[Theorem~1.1(i)]{PopaRusJvN}, stated there for real functions. The closed form (\ref{eq:Q2}) is not part of that statement.

We close the list of examples with a weighted family. Let $R$ be a commutative ring and let $H$ be an $R$-module. For maps $\mu\colon G^2\to G$ and $\omega\colon G^2\to R$, the operator $\omega\mu^{*}\colon \mathcal{H}_1\to \mathcal{H}_2$ defined by
\[
(\omega\mu^{*})[f](x,y) := \omega(x,y)\,f\bigl(\mu(x,y)\bigr)
\]
is called the \emph{weighted composition operator} of $\mu$ with weight~$\omega$. For $\mu$ the addition of a semigroup, we call $\omega\mathsf{E}$ a \emph{twisted translation}, after the twisted group ring of Remark~\ref{rem:twisted}. Proposition~\ref{prop:twisted} gives a sufficient condition for the associativity of $\omega\mu^{*}$ over a general~$R$. For $R=\mathbb{Z}$ and a target that no non-zero integer annihilates, it also gives two partial converses. Section~\ref{sec:image} uses a general~$R$.

\begin{proposition}[weighted composition operators]\label{prop:twisted}
Let $\mu\colon G^2\to G$, let $R$ be a commutative ring, let $H$ be an $R$-module and let $\omega\colon G^2\to R$.
\begin{enumerate}
\item[\rm (i)] If $\mu$ is associative and, for all $x,y,z\in G$,
\begin{equation}\label{eq:mult-cocycle}
\omega(x,y)\,\omega\bigl(\mu(x,y),z\bigr) = \omega(y,z)\,\omega\bigl(x,\mu(y,z)\bigr),
\end{equation}
then $\omega\mu^{*}$ is associative.
\item[\rm (ii)] Let $R=\mathbb{Z}$ and let $H$ be an abelian group of infinite exponent, that is, one that no non-zero integer annihilates. If $\mu$ is associative and $\omega\mu^{*}$ is associative, then $\omega$ satisfies (\ref{eq:mult-cocycle}). If $\omega$ vanishes nowhere and $\omega\mu^{*}$ is associative, then $\mu$ is associative.
\end{enumerate}
\end{proposition}

\begin{proof}
From (\ref{eq:lift}), the lift of $\omega\mu^{*}$ in a slot folds the active pair by $\mu$ and multiplies by its weight, so
\[
\begin{gathered}
(\omega\mu^{*})^{[1]}(\omega\mu^{*})[f](x,y,z) = \omega(x,y)\,\omega\bigl(\mu(x,y),z\bigr)\, f\bigl(\mu(\mu(x,y),z)\bigr),\\
(\omega\mu^{*})^{[2]}(\omega\mu^{*})[f](x,y,z) = \omega(y,z)\,\omega\bigl(x,\mu(y,z)\bigr)\, f\bigl(\mu(x,\mu(y,z))\bigr).
\end{gathered}
\]
Fix $(x,y,z)$ and write $\omega_1$, $\omega_2$ for the two weights and $g_1$, $g_2$ for the two arguments of~$f$ in these formulas.

\emph{Part {\rm(i)}.} If $\mu$ is associative and $\omega$ satisfies (\ref{eq:mult-cocycle}), then $g_1=g_2$ and $\omega_1=\omega_2$ at every triple, and the two iterates coincide.

\emph{Part {\rm(ii)}, first assertion.} Let $\mu$ and $\omega\mu^{*}$ be associative. Then $g_1=g_2$, and the two iterates agree at $(x,y,z)$ for every $f$ exactly when $(\omega_1-\omega_2)c=0$ for every $c\in H$, as one sees from the function with value $c$ at $g_1$ and $0$ elsewhere. Since no non-zero integer annihilates $H$, this forces $\omega_1=\omega_2$, which is (\ref{eq:mult-cocycle}).

\emph{Part {\rm(ii)}, second assertion.} If $g_1\neq g_2$ at some triple, the function with value $c$ at $g_1$ and $0$ elsewhere gives $\omega_1 c=0$ for every $c\in H$, hence $\omega_1=0$ by the hypothesis on $H$, contradicting the hypothesis on~$\omega$. So $g_1=g_2$ at every triple.
\end{proof}

\begin{remark}\label{rem:twisted}
Identity (\ref{eq:mult-cocycle}) is the associativity condition of a twisted semigroup ring. For $\mu$ the addition of an abelian group it is the classical associativity condition of a twisted group ring \cite[Ch.~1, §2]{Passman1977}, the multiplicative counterpart of the additive cocycles of \cite{EilenbergMacLane}. Proposition~\ref{prop:twisted} states these conditions for operators on the hierarchy.
\end{remark}

\section{The tower theorem}\label{sec:towers}

In Proposition~\ref{prop:R} and Theorem~\ref{thm:towers}, $G$ and $H$ are arbitrary non-empty sets and $\mathsf{T}\colon \mathcal{H}_1\to \mathcal{H}_2$ is an associative operator. The other results of this section state their own hypotheses.

\begin{proposition}[slot-wise propagation]\label{prop:R}
Let $\mathsf{T}$ be associative, and let $n\ge 1$ and $1\le i\le n$. Then, as maps $\mathcal{H}_n\to \mathcal{H}_{n+2}$,
\begin{equation}\label{eq:R}
\mathsf{T}^{[i]}\,\mathsf{T}^{[i]} \;=\; \mathsf{T}^{[i+1]}\,\mathsf{T}^{[i]}.
\end{equation}
\end{proposition}

\begin{proof}
Let $\psi\in \mathcal{H}_n$ and $(x_1,\dots,x_{n+2})\in G^{n+2}$. Define $f\in \mathcal{H}_1$ by
\[
f(x) := \psi(x_1,\dots,x_{i-1},x,x_{i+3},\dots,x_{n+2}).
\]
From (\ref{eq:lift}), applied twice, the two sides of~(\ref{eq:R}) evaluated at $(x_1,\dots,x_{n+2})$ become
\[
\mathsf{T}\bigl[\mathsf{T}[f](\,\cdot\,,x_{i+2})\bigr](x_i,x_{i+1}) = \bigl(\mathsf{T}^{[1]}\mathsf{T}\bigr)[f](x_i,x_{i+1},x_{i+2})
\]
and
\[
\mathsf{T}\bigl[\mathsf{T}[f](x_i,\,\cdot\,)\bigr](x_{i+1},x_{i+2}) = \bigl(\mathsf{T}^{[2]}\mathsf{T}\bigr)[f](x_i,x_{i+1},x_{i+2}),
\]
and these coincide by associativity.
\end{proof}

\begin{definition}
Let $\mathsf{T}\colon \mathcal{H}_1\to \mathcal{H}_2$ be an operator, let $n\ge 1$, and let $j_1,\dots,j_n$ be integers with $1\le j_k\le k$ for every~$k$. The composite
\[
\mathsf{T}^{[j_{n}]}\cdots \mathsf{T}^{[j_2]}\,\mathsf{T}^{[j_1]} \colon \mathcal{H}_1\longrightarrow \mathcal{H}_{n+1}
\]
is called a \emph{tower of length $n$} over $\mathsf{T}$.
\end{definition}

\begin{remark}
The factor applied first, $\mathsf{T}^{[j_1]}$, acts on $\mathcal{H}_1$, which has a single slot, so necessarily $j_1=1$ and $\mathsf{T}^{[j_1]}=\mathsf{T}$. The only tower of length one is $\mathsf{T}$ itself.
\end{remark}

\begin{theorem}[tower theorem]\label{thm:towers}
Let $\mathsf{T}$ be associative. Then all towers of the same length coincide. The common value of the towers of length~$n$ is denoted $\mathsf{T}^{\,n}\colon \mathcal{H}_1\to \mathcal{H}_{n+1}$, and for every $i$ with $1\le i\le n+1$,
\begin{equation}\label{eq:power-rec}
\mathsf{T}^{\,n+1} = \mathsf{T}^{[i]}\,\mathsf{T}^{\,n}.
\end{equation}
\end{theorem}

\begin{proof}
We prove by induction on $n\ge 1$ the statement: all towers of length~$n$ coincide, and their common value $\mathsf{T}^{\,n}$ satisfies (\ref{eq:power-rec}) one level up, that is, $\mathsf{T}^{[i]}\mathsf{T}^{\,n}$ is independent of $i\in\{1,\dots,n+1\}$.

For $n=1$ the only tower is $\mathsf{T}$, and the independence of $\mathsf{T}^{[i]}\mathsf{T}$ from $i\in\{1,2\}$ is the definition of associativity.

Assume the statement for $n-1$, so that all towers of length $n-1$ equal $\mathsf{T}^{\,n-1}$ and $\mathsf{T}^{\,n} := \mathsf{T}^{[i]}\mathsf{T}^{\,n-1}$ does not depend on~$i$. Every tower of length $n$ is of the form $\mathsf{T}^{[j]}\mathsf{S}$ with $\mathsf{S}$ a tower of length $n-1$, hence equals $\mathsf{T}^{[j]}\mathsf{T}^{\,n-1} = \mathsf{T}^{\,n}$. All towers of length $n$ therefore coincide. It remains to prove that $\mathsf{T}^{[i]}\mathsf{T}^{\,n}$ is independent of~$i$, and it suffices to compare consecutive indices. Fix $1\le i\le n$. By the induction hypothesis, $\mathsf{T}^{\,n} = \mathsf{T}^{[i]}\mathsf{T}^{\,n-1}$, so, by Proposition~\ref{prop:R} applied on $\mathcal{H}_n$,
\[
\mathsf{T}^{[i]}\,\mathsf{T}^{\,n}
= \mathsf{T}^{[i]}\,\mathsf{T}^{[i]}\,\mathsf{T}^{\,n-1}
= \mathsf{T}^{[i+1]}\,\mathsf{T}^{[i]}\,\mathsf{T}^{\,n-1}
= \mathsf{T}^{[i+1]}\,\mathsf{T}^{\,n}. \qedhere
\]
\end{proof}

\begin{remark}
Theorem~\ref{thm:towers} is the operator analogue of the generalised associative law, by which all bracketings of a product of $n+1$ factors under an associative operation coincide \cite[Ch.~I, §1, no.~3]{Bourbaki}.
\end{remark}

\begin{remark}\label{rem:no-exchange}
Theorem~\ref{thm:towers} uses only the relation (\ref{eq:R}), which involves a slot and the next one. For lifts in two slots $i<j$ there is a second relation, which an operator may or may not satisfy, the \emph{exchange law}
\[
\mathsf{T}^{[i]}\,\mathsf{T}^{[j]} = \mathsf{T}^{[j+1]}\,\mathsf{T}^{[i]},\qquad 1\le i<j\le n,
\]
an identity between maps $\mathcal{H}_n\to\mathcal{H}_{n+2}$. The relation (\ref{eq:R}) and the exchange law together are the relations satisfied by the inner coface maps of a cosimplicial object \cite[VII.5]{MacLaneCWM}. The composition operators satisfy the exchange law, whether or not their operations are associative. For any two maps $\mu,\nu\colon G^2\to G$ and any $1\le i<j\le n$, the maps $\mu_{[i]}\circ\nu_{[j+1]}$ and $\nu_{[j]}\circ\mu_{[i]}$ from $G^{n+2}$ to $G^{n}$ coincide, because both apply $\mu$ to the slots $i,i+1$ and $\nu$ to the slots $j+1,j+2$ of the argument. By Lemma~\ref{lem:basic}(i) this gives $(\mu^{*})^{[i]}(\nu^{*})^{[j]}=(\nu^{*})^{[j+1]}(\mu^{*})^{[i]}$, and in particular the exchange law for every composition operator. When $H$ is an abelian group, the exchange law for an integer combination $\sum_k m_k\mu_k^{*}$ follows by expanding both sides with Lemma~\ref{lem:basic}(ii) and~(iii) and comparing them term by term. The same argument applies to the weighted composition operators, since the two weights are multiplied in either order. Associativity alone does not imply the exchange law, as Example~\ref{ex:no-exchange} shows.
\end{remark}

\begin{example}[associativity does not imply the exchange law]\label{ex:no-exchange}
Let $\beta\colon H^2\to H$ be a map and let $\mathsf{T}[u](x,y):=\beta(u(x),u(y))$. The operator $\mathsf{T}$ is associative as soon as, for all $a,b,c\in H$,
\[
\beta\bigl(\beta(a,c),\beta(b,c)\bigr) = \beta\bigl(\beta(a,b),\beta(a,c)\bigr),
\]
because for $a=u(x)$, $b=u(y)$ and $c=u(z)$ the two sides are the values of $\mathsf{T}^{[1]}\mathsf{T}[u]$ and $\mathsf{T}^{[2]}\mathsf{T}[u]$ at $(x,y,z)$. If $G$ has at least two elements, the exchange law $\mathsf{T}^{[1]}\mathsf{T}^{[2]}=\mathsf{T}^{[3]}\mathsf{T}^{[1]}$ on $\mathcal{H}_2$ holds if and only if $\beta$ is medial, that is, $\beta(\beta(a,b),\beta(c,d)) = \beta(\beta(a,c),\beta(b,d))$ for all $a,b,c,d\in H$. On the three-element set $H=\{0,1,2\}$, let $\beta$ agree with the first projection except for $\beta(0,2)=2$. This operation satisfies the first identity and fails the second at $(a,b,c,d)=(0,0,1,2)$, where the two sides are $0$ and~$2$. So $\mathsf{T}$ is associative and, when $G$ has at least two elements, it does not obey the exchange law. The action of an associative operator on the hierarchy need not be cosimplicial. A sequel in preparation delimits the operators for which the exchange law holds.
\end{example}

\begin{example}\label{ex:powers}
The powers of the composition operators of Example~\ref{ex:comp} are computed by induction from Lemma~\ref{lem:basic}(i). For an associative operation $\mu$, in the notation of Section~\ref{sec:lifts},
\[
(\mu^{*})^{\,n}[f](x_1,\dots,x_{n+1}) = f\bigl(\mu(x_1,\dots,x_{n+1})\bigr).
\]
Thus, over a semigroup, $\mathsf{E}^{\,n}[f] = f(x_1+\dots+x_{n+1})$. Over any set, $\mathsf{J}_1^{\,n}[f] = f(x_1)$ and $\mathsf{J}_2^{\,n}[f] = f(x_{n+1})$. Over a pointed set, $\mathsf{J}_e^{\,n}[f] = f(e)$. Over a linearly ordered set, $(\max{}^{*})^{\,n}[f] = f(\max(x_1,\dots,x_{n+1}))$. In each case the argument is $(x_1,\dots,x_{n+1})$. For a weighted composition operator with $\mu$ associative and $\omega$ satisfying (\ref{eq:mult-cocycle}), which is associative by Proposition~\ref{prop:twisted}(i), induction along (\ref{eq:power-rec}) gives
\[
(\omega\mu^{*})^{\,n}[f](x_1,\dots,x_{n+1}) = \Bigl(\,\prod_{k=1}^{n}\omega\bigl(\mu(x_1,\dots,x_k),\,x_{k+1}\bigr)\Bigr)\, f\bigl(\mu(x_1,\dots,x_{n+1})\bigr).
\]
By Theorem~\ref{thm:towers} the left-hand side does not depend on the order in which the lifts are taken. The different orders lead to the products of weights along the different bracketings of $x_1,\dots,x_{n+1}$, and these products coincide by repeated use of (\ref{eq:mult-cocycle}). This is the generalised associative law in the twisted semigroup ring. For the Cauchy difference, over a semigroup $G$ and an abelian group $H$, induction along (\ref{eq:power-rec}) gives
\[
\mathsf{L}^{\,n}[f](x_1,\dots,x_{n+1}) = \sum_{\emptyset\neq S\subseteq\{1,\dots,n+1\}} (-1)^{n+1-|S|}\, f\Bigl(\sum_{k\in S} x_k\Bigr),
\]
the \emph{Cauchy difference of order} $n$, in the terminology of \cite{GuoMaLi,Ng2015,LiNg2016}, whose kernel is described for $n=2$ in \cite{NgZhao2013,Stetkaer2017} and, at higher orders over particular classes of groups and semigroups, in \cite{GuoMaLi,Ng2015,LiNg2016}. The inner sum runs over $S$ in increasing order of the indices, so the formula holds over any semigroup. The recursion by which \cite{GuoMaLi} defines the difference of order $n$ from the difference of order $n-1$ is the lift $\mathsf{L}^{[1]}$, so that the tower of $\mathsf{L}$ lifted always in the first slot is implicit in that definition, and Theorem~\ref{thm:towers} adds that the slot of the lift is immaterial. The case $n=2$ is Proposition~\ref{prop:L}. We keep the name \emph{Fréchet operator}, and the letter~$\mathsf{F}$, for the second-order difference $\mathsf{F}=\mathsf{L}^{\,2}$ alone. Remark~\ref{rem:two-towers} relates it to Fréchet's functional equation \cite{Frechet1909}.
\end{example}

\begin{example}\label{ex:PQ-powers}
Let $G$ and $H$ be abelian groups. For the d'Alembert operator, induction along (\ref{eq:power-rec}) with $i=n+1$ gives
\[
\mathsf{P}^{\,n}[f](x_1,\dots,x_{n+1}) = \sum_{\varepsilon\in\{\pm 1\}^{n}} f(x_1+\varepsilon_1 x_2+\dots+\varepsilon_n x_{n+1}).
\]
Indeed, the lift in the last slot replaces the variable $x_{n+1}$ by $x_{n+1}\pm x_{n+2}$, which doubles the number of sign vectors. For the Jordan--von Neumann operator, Proposition~\ref{prop:Q} and Theorem~\ref{thm:towers} give the powers $\mathsf{Q}^{\,n}\colon \mathcal{H}_1\to \mathcal{H}_{n+1}$, of which $\mathsf{Q}^{\,2}$ is given by~(\ref{eq:Q2}). Proposition~\ref{prop:Qn} computes them in closed form.
\end{example}

\begin{proposition}\label{prop:Qn}
Let $G$ and $H$ be abelian groups. For all $n\ge 1$ and $f\in \mathcal{H}_1$,
\begin{equation}\label{eq:Qn}
\mathsf{Q}^{\,n}[f](x_1,\dots,x_{n+1}) \;=\; \sum_{\emptyset\neq S\subseteq\{1,\dots,n+1\}} (-2)^{\,n+1-|S|} \sum_{\varepsilon\in\{\pm1\}^{S\setminus\{m_S\}}} f\Bigl(x_{m_S}+\sum_{k\in S\setminus\{m_S\}}\varepsilon_k\, x_k\Bigr),
\end{equation}
where $m_S := \min S$.
\end{proposition}

\begin{proof}
We argue by induction on $n$ along (\ref{eq:power-rec}) with $i=n+1$, that is, lifting in the last slot. The base $n=1$ is the definition of~$\mathsf{Q}$. Assume (\ref{eq:Qn}) at level $n$ and write $\psi := \mathsf{Q}^{\,n}[f]$. By (\ref{eq:lift}),
\begin{multline*}
\mathsf{Q}^{\,n+1}[f](x_1,\dots,x_{n+2}) = \psi(x_1,\dots,x_n,x_{n+1}+x_{n+2}) + \psi(x_1,\dots,x_n,x_{n+1}-x_{n+2})\\
- 2\psi(x_1,\dots,x_n,x_{n+1}) - 2\psi(x_1,\dots,x_n,x_{n+2}).
\end{multline*}
Each of the four evaluations substitutes, in every term of (\ref{eq:Qn}), the variable $x_{n+1}$ by $x_{n+1}+x_{n+2}$, by $x_{n+1}-x_{n+2}$, by $x_{n+1}$ and by $x_{n+2}$ respectively, with the coefficients $1$, $1$, $-2$, $-2$. We call these the four substitutions. The terms of (\ref{eq:Qn}) at level $n$ are indexed by pairs $(S,\varepsilon)$ with $\emptyset\neq S\subseteq\{1,\dots,n+1\}$, those at level $n+1$ by pairs $(S',\varepsilon')$ with $\emptyset\neq S'\subseteq\{1,\dots,n+2\}$, and we sort the terms produced by the four substitutions according to $S'\cap\{n+1,n+2\}$. A term of $\psi$ with $n+1\notin S$ does not involve $x_{n+1}$ and is left unchanged by every substitution. A term with $n+1\in S$ produces, under the four substitutions in turn, terms with both $n+1$ and $n+2$, both again, only $n+1$, and only $n+2$ in the index set. Four cases result.

\emph{Neither $n+1$ nor $n+2$ lies in $S'$.} Such a term comes only from $S=S'$, with $n+1\notin S$, under all four substitutions, which leave it unchanged. Its total coefficient is $(1+1-2-2)\,(-2)^{\,n+1-|S|} = (-2)^{\,n+2-|S'|}$, and its signs are those of $(S,\varepsilon)$, as required.

\emph{Both $n+1$ and $n+2$ lie in $S'$.} Such a term comes only from $S=S'\setminus\{n+2\}$, which contains $n+1$, under the first two substitutions. In the term of $\psi$ the variable $x_{n+1}$ carries the sign $\varepsilon_{n+1}$, or the sign $+$ when $S=\{n+1\}$. After the substitution $x_{n+1}\mapsto x_{n+1}\pm x_{n+2}$ the variable $x_{n+2}$ carries the sign $\pm\varepsilon_{n+1}$, so that, as the two substitutions are added, this sign runs over both values. The two substitutions together produce exactly the sum over $\varepsilon'_{n+2}\in\{\pm1\}$ required at level $n+1$, with $m_{S'}=m_S$, the remaining signs unchanged, and the coefficient $(-2)^{\,n+1-|S|} = (-2)^{\,n+2-|S'|}$, since $|S'|=|S|+1$.

\emph{Only $n+1$ lies in $S'$.} Such a term comes only from $S=S'$, which contains $n+1$, under the third substitution, which leaves the term unchanged and multiplies it by $-2$. Its coefficient is $-2\,(-2)^{\,n+1-|S|} = (-2)^{\,n+2-|S'|}$, with the signs unchanged.

\emph{Only $n+2$ lies in $S'$.} Such a term comes only from $S=(S'\setminus\{n+2\})\cup\{n+1\}$ under the fourth substitution, $x_{n+1}\mapsto x_{n+2}$, which transfers the sign of $x_{n+1}$ to $x_{n+2}$. When $S=\{n+1\}$ it turns $f(x_{n+1})$ into $f(x_{n+2})$, in agreement with $m_{S'}=n+2$ in that case. Its coefficient is $-2\,(-2)^{\,n+1-|S|} = (-2)^{\,n+2-|S'|}$, since $|S'|=|S|$.

In each case the terms produced are exactly the terms of (\ref{eq:Qn}) at level $n+1$ indexed by the pairs $(S',\varepsilon')$ of that case, each reached once. This completes the induction.
\end{proof}

\begin{remark}
Formula (\ref{eq:Qn}) is related to the formulas of Examples \ref{ex:powers} and~\ref{ex:PQ-powers}. Its term $S=\{1,\dots,n+1\}$ is the d'Alembert sum $\mathsf{P}^{\,n}[f]$. Replacing the coefficients $(-2)^{\,n+1-|S|}$ by $(-1)^{\,n+1-|S|}$, and each inner sum by its term with all signs positive, gives the Cauchy difference $\mathsf{L}^{\,n}[f]$.
\end{remark}

\subsection*{The towers of \texorpdfstring{$\Delta$}{Delta} and the calculus at the origin}
For an operator that is not associative, the towers of a given length need not coincide. For $\Delta$ all of them can be described, in terms of the classical difference calculus. Let $G$ be a commutative monoid and $H$ an abelian group. For $t\in G$, the \emph{shift} and the \emph{difference with increment}~$t$ are the operators $\mathsf{E}_t,\Delta_t\colon \mathcal{H}_1\to \mathcal{H}_1$ of the finite calculus \cite{Jordan},
\begin{equation}\label{eq:shift-diff}
\mathsf{E}_t[f](x) := f(t+x), \qquad \Delta_t[f](x) := f(t+x)-f(x),
\end{equation}
so that $\mathsf{E}[f](x,y) = \mathsf{E}_x[f](y)$ and $\Delta[f](x,y) = \Delta_x[f](y)$. Thus the first variable of $\Delta[f]$ is the increment. The family $(\Delta_t)_{t\in G}$ is commutative, $\Delta_s\Delta_t=\Delta_t\Delta_s$, and it commutes with the shifts, $\mathsf{E}_s\Delta_t=\Delta_t\mathsf{E}_s$. It satisfies $\Delta_0=0$ and $\mathsf{E}_t=\mathsf{I}+\Delta_t$, where $\mathsf{I}$ is the identity operator of $\mathcal{H}_1$. For all $s,t\in G$ it obeys the addition formula
\begin{equation}\label{eq:addition}
\Delta_{s+t} = \Delta_s + \Delta_t + \Delta_s\Delta_t,
\end{equation}
since both sides send $f$ to the function $x\mapsto f(s+t+x)-f(x)$, as a direct expansion shows.

Let $G$ be a commutative monoid, let $H$ be an abelian group and let $n\ge 1$. For a subset $A\subseteq\{2,\dots,n\}$ put $B:=\{2,\dots,n+1\}\setminus A$. The operator $\mathsf{D}_A\colon \mathcal{H}_1\to \mathcal{H}_{n+1}$ defined by
\[
\mathsf{D}_A[f](x_1,\dots,x_{n+1}) := \Bigl(\Delta_{x_1}\prod_{k\in A}\Delta_{x_k}[f]\Bigr)\Bigl(\sum_{k\in B}x_k\Bigr)
\]
is called a \emph{partial mixed difference}. Its increments are $x_1$ and the $x_k$ with $k\in A$, and it is evaluated at the sum of the remaining variables. The arity $n+1$ is not recorded in the notation and is supplied by the context.

\begin{proposition}[towers of $\Delta$]\label{prop:delta-towers}
Let $G$ be a commutative monoid, let $H$ be an abelian group and let $n\ge 1$. The towers of $\Delta$ of length $n$ are exactly the operators $\mathsf{D}_A$, $A\subseteq\{2,\dots,n\}$. The tower lifted always in the first slot is $\mathsf{D}_{\emptyset}$, the first difference $f({x_1+\dots+x_{n+1}})-f({x_2+\dots+x_{n+1}})$. The tower $\Delta^{[n]}\cdots\Delta^{[2]}\Delta$, lifted always in the last slot, is $\mathsf{D}_{\{2,\dots,n\}}$, the full mixed difference
\begin{equation}\label{eq:delta-mixed}
\mathsf{D}_{\{2,\dots,n\}}[f](x_1,\dots,x_{n+1}) = (\Delta_{x_1}\cdots\Delta_{x_n}[f])(x_{n+1}).
\end{equation}
If, moreover, $G$ contains an element $x_0$ whose multiples $kx_0$, $0\le k\le n-1$, are pairwise distinct (for instance an element of infinite order) and $H$ is non-trivial, the operators $\mathsf{D}_A$ are pairwise distinct, so that $\Delta$ has exactly $2^{\,n-1}$ towers of length~$n$.
\end{proposition}

\begin{proof}
For $n=1$ the only tower is $\Delta$, and
\[
\Delta[f](x_1,x_2) = f(x_1+x_2)-f(x_2) = (\Delta_{x_1}[f])(x_2) = \mathsf{D}_{\emptyset}[f](x_1,x_2).
\]
Assume that every tower of length $n$ is one of the operators $\mathsf{D}_A$, $A\subseteq\{2,\dots,n\}$, and lift $\mathsf{D}_A$ in a slot $j\in\{1,\dots,n+1\}$. By (\ref{eq:lift}),
\begin{multline*}
\Delta^{[j]}[\psi](x_1,\dots,x_{n+2}) = \psi(x_1,\dots,x_{j-1},\,x_j+x_{j+1},\,x_{j+2},\dots,x_{n+2})\\
- \psi(x_1,\dots,x_{j-1},\,x_{j+1},\,x_{j+2},\dots,x_{n+2}),
\end{multline*}
so that the variables of $\psi$ after the active pair are shifted by one place. Let $\psi=\mathsf{D}_A[f]$. If the $j$-th variable of $\psi$ is an increment, the two values differ only in that increment, and (\ref{eq:addition}) gives $\Delta_{x_j+x_{j+1}}-\Delta_{x_{j+1}} = \Delta_{x_j}+\Delta_{x_j}\Delta_{x_{j+1}} = \Delta_{x_j}\mathsf{E}_{x_{j+1}}$. Since the shifts commute with the differences, $(\Delta_{x_j}\mathsf{E}_{x_{j+1}}[u])(s) = (\Delta_{x_j}[u])(x_{j+1}+s)$ for every $u\in \mathcal{H}_1$ and every $s\in G$, so the increments are unchanged and the new variable $x_{j+1}$ joins the evaluation sum. If the $j$-th variable of $\psi$ belongs to the evaluation sum, which is then $x_j+s'$ with $s'$ the sum of its other terms, the two values are $u({x_j+x_{j+1}+s'})-u({x_{j+1}+s'}) = (\Delta_{x_j}[u])(x_{j+1}+s')$, where $u$ is the mixed difference of $f$ with the increments of $\psi$: the variable $x_j$ becomes an increment and $x_{j+1}$ joins the sum. In both cases the result is $\mathsf{D}_{A'}[f]$ for a subset $A'\subseteq\{2,\dots,n+1\}$, because $x_1$ remains an increment (it is untouched when $j\ge 2$, and it stays an increment when $j=1$) and the last variable $x_{n+2}$ belongs to the sum (it is the shifted last variable of $\psi$ when $j\le n$, and the new variable $x_{j+1}$ when $j=n+1$). Hence every tower of length $n+1$ is one of the $\mathsf{D}_{A'}$. The two cases also identify the extreme towers: lifting always in the first slot, an increment slot, keeps $A=\emptyset$, and lifting always in the last slot, which belongs to the sum, produces $A=\{2,\dots,n\}$, which is (\ref{eq:delta-mixed}).

Conversely, every $\mathsf{D}_A$ is a tower, by induction on $n$. For $A=\emptyset$ it is the tower lifted always in the first slot, an increment slot, which keeps $A=\emptyset$ at every length. If $A\neq\emptyset$, let $\ell:=\max A$. By induction $\mathsf{D}_{A\setminus\{\ell\}}$, with $A\setminus\{\ell\}\subseteq\{2,\dots,\ell-1\}$, is a tower of length $\ell-1$, at which length the variable $x_{\ell}$ belongs to the sum. Lifting in slot $\ell$ turns $x_{\ell}$ into an increment and produces $\mathsf{D}_A$ at length $\ell$, and lifting $n-\ell$ further times in slot $\ell$, now an increment slot, keeps $A$ and produces $\mathsf{D}_A$ at length $n$.

Finally, let $A\neq A'$. Interchanging them if necessary, pick $i\in A\setminus A'$, and let $x_0$ be as in the statement. At every point with $x_i=0$ the operator $\mathsf{D}_A$ vanishes, because the increments commute and $\Delta_0=0$. Consider the point with $x_1=x_0$, $x_k=x_0$ for $k\in A'$, and $x_k=0$ for all other $k$, including $k=i$, which lies in $\{2,\dots,n+1\}\setminus A'$. There $\mathsf{D}_{A'}[f] = (\Delta_{x_0}^{\,m}[f])(0) = \sum_{k=0}^{m}(-1)^{m-k}\binom{m}{k}f(kx_0)$ with $m=|A'|+1\le n-1$, because $A'$ is contained in $\{2,\dots,n\}$ and does not contain~$i$, and for the function $f$ with value $c\neq 0$ at $mx_0$ and $0$ elsewhere this equals $c$, because the points $0,x_0,\dots,mx_0$ are pairwise distinct. Thus $\mathsf{D}_A\neq\mathsf{D}_{A'}$, and the $2^{\,n-1}$ subsets $A$ give $2^{\,n-1}$ distinct towers.
\end{proof}

\begin{remark}[the hypothesis on $x_0$]\label{rem:delta-sharp}
The hypothesis on $x_0$ in Proposition~\ref{prop:delta-towers} cannot be dropped. Over $G=H=\mathbb{Z}_2$ every $\mathsf{D}_A$ with $A\neq\emptyset$ vanishes identically, since $\Delta_1\Delta_1[f](s) = 2f(s)-2f(s+1) = 0$ for every~$f$. The towers are then $\mathsf{D}_{\emptyset}$, which does not vanish, and the zero operator, so $\Delta$ has exactly two distinct towers of each length $n\ge 2$, and the count $2^{\,n-1}$ fails from $n=3$ on. On the other hand, the hypothesis is not necessary. Over $G=\mathbb{Z}_2$ and $H=\mathbb{Z}$ one has $\Delta_1\Delta_1=-2\Delta_1$, hence $(\Delta_1^{\,m}[f])(0)=(-2)^{\,m-1}\bigl(f(1)-f(0)\bigr)$ for $m\ge 1$. The last step of the proof of Proposition~\ref{prop:delta-towers}, applied with $x_0=1$ and with a function $f$ such that $f(0)=0$ and $f(1)=1$, then shows that the operators $\mathsf{D}_A$ are again pairwise distinct.
\end{remark}

Evaluating the mixed differences at the origin gives a sequence of operators which, unlike the towers of $\Delta$, are, from the second term on, the powers of a single associative operator. For a commutative monoid $G$, an abelian group $H$ and $n\ge 1$, the operator $\mathsf{D}_{\langle n\rangle}\colon \mathcal{H}_1\to \mathcal{H}_n$ defined by
\begin{equation}\label{eq:Dn}
\mathsf{D}_{\langle n\rangle}[f](x_1,\dots,x_n) := (\Delta_{x_1}\cdots\Delta_{x_n}[f])(0)
= \sum_{S\subseteq\{1,\dots,n\}} (-1)^{n-|S|}\, f\Bigl(\sum_{k\in S} x_k\Bigr)
\end{equation}
is called the \emph{mixed difference at the origin}. The expansion follows by induction on~$n$ from (\ref{eq:shift-diff}), with the convention that an empty sum of elements of $G$ is $0$, so that the term $S=\emptyset$ is $(-1)^{n}f(0)$. Since $\Delta_0=0$ and the increments commute, $\mathsf{D}_{\langle n\rangle}[f]$ vanishes whenever one of its arguments is~$0$. These operators are the left-hand sides of Fréchet's functional equation \cite{Frechet1909}, whose homogeneous theory, the description of the functions with vanishing differences of a given order, is developed in \cite{MazurOrlicz,Djokovic,Szekelyhidi2014}. For the inhomogeneous equation $\mathsf{D}_{\langle n\rangle}[f]=\psi$, Example~\ref{ex:coherent} gives necessary conditions on the datum~$\psi$. For $n=2$ the definition reads $\mathsf{D}_{\langle 2\rangle}[f](x,y) = f(x+y)-f(x)-f(y)+f(0)$, that is,
\begin{equation*}
\mathsf{D}_{\langle 2\rangle} = \mathsf{L} + \mathsf{J}_0.
\end{equation*}

\begin{proposition}\label{prop:D2}
Let $G$ be a monoid and $H$ an abelian group. Every integer combination $m\mathsf{L}+n\mathsf{J}_0$ is associative. In particular the operator $\mathsf{L}+\mathsf{J}_0$ is associative, and if $G$ is moreover commutative it is the mixed difference $\mathsf{D}_{\langle 2\rangle}$, with $\mathsf{D}_{\langle 2\rangle}^{\,2} = \mathsf{D}_{\langle 3\rangle}$.
\end{proposition}

\begin{proof}
Write $\mathsf{T}=m\mathsf{L}+n\mathsf{J}_0$ and $\chi=\mathsf{T}[f]$, so that $\chi(x,y)=m\bigl[f(x+y)-f(x)-f(y)\bigr]+nf(0)$ for all $x,y\in G$, and $\chi(0,y)=\chi(x,0)=(n-m)f(0)$. By Lemma~\ref{lem:basic}(ii),
\[
\begin{gathered}
\mathsf{T}^{[1]}[\chi](x,y,z) = m\chi(x{+}y,z)-m\chi(x,z)-m\chi(y,z)+n\chi(0,z),\\
\mathsf{T}^{[2]}[\chi](x,y,z) = m\chi(x,y{+}z)-m\chi(x,y)-m\chi(x,z)+n\chi(x,0).
\end{gathered}
\]
The terms $-m\chi(x,z)$ and the constant terms agree, and both
\[
\chi(x{+}y,z)-\chi(y,z) \qquad\text{and}\qquad \chi(x,y{+}z)-\chi(x,y)
\]
expand to $m\bigl[f(x+y+z)-f(x+y)-f(y+z)+f(y)\bigr]$, so the two iterates coincide.

For the second assertion, write $\psi=\mathsf{D}_{\langle 2\rangle}[f]$, so that $\psi(x,y)=f(x+y)-f(x)-f(y)+f(0)$ for all $x,y\in G$, a function that vanishes when either argument is $0$. Expanding the first iterate as in the proof of Proposition~\ref{prop:L}, with the extra term $\psi(0,z)=0$ contributed by the summand $\mathsf{J}_0$,
\[
\mathsf{D}_{\langle 2\rangle}^{[1]}\mathsf{D}_{\langle 2\rangle}[f](x,y,z) = \psi(x+y,z)-\psi(x,z)-\psi(y,z)+\psi(0,z).
\]
Substituting that formula and simplifying, this expands to
\[
f(x{+}y{+}z)-f(x{+}y)-f(y{+}z)-f(x{+}z)+f(x)+f(y)+f(z)-f(0) = \mathsf{D}_{\langle 3\rangle}[f](x,y,z),
\]
by (\ref{eq:Dn}) at $n=3$. The second iterate is $\psi(x,y+z)-\psi(x,y)-\psi(x,z)+\psi(x,0)$, and the same substitution expands it to the same eight-term sum. The two iterates coincide, and their common value is $\mathsf{D}_{\langle 3\rangle}[f]$.
\end{proof}

The powers of $\mathsf{D}_{\langle 2\rangle}$ are obtained from a substitution rule, which follows from the addition formula (\ref{eq:addition}).

\begin{lemma}[substitution rule]\label{lem:substitution}
Let $G$ be a commutative monoid and $H$ an abelian group. For all $n\ge 1$ and $1\le i\le n$,
\[
\mathsf{L}^{[i]}\,\mathsf{D}_{\langle n\rangle} = \mathsf{D}_{\langle n+1\rangle}, \qquad \mathsf{J}_0^{[i]}\,\mathsf{D}_{\langle n\rangle} = 0.
\]
\end{lemma}

\begin{proof}
By (\ref{eq:lift}) and (\ref{eq:Dn}), the three evaluations of $\mathsf{D}_{\langle n\rangle}[f]$ that $\mathsf{L}^{[i]}$ combines differ only in one increment, so, since the increments commute,
\[
\mathsf{L}^{[i]}\,\mathsf{D}_{\langle n\rangle}[f](x_1,\dots,x_{n+1})
= \Bigl(\bigl(\Delta_{x_i+x_{i+1}}-\Delta_{x_i}-\Delta_{x_{i+1}}\bigr)\prod_{k\notin\{i,i+1\}}\Delta_{x_k}\,[f]\Bigr)(0),
\]
the product running over the remaining increments, of which there are none when $n=1$. By (\ref{eq:addition}) the bracket equals $\Delta_{x_i}\Delta_{x_{i+1}}$, so the right-hand side is $\mathsf{D}_{\langle n+1\rangle}[f](x_1,\dots,x_{n+1})$. For the second identity, $\mathsf{J}_0^{[i]}$ substitutes $0$ into the $i$-th slot, by Lemma~\ref{lem:basic}(i), and $\mathsf{D}_{\langle n\rangle}[f]$ vanishes there.
\end{proof}

\begin{proposition}[calculus at the origin as a tower]\label{prop:origin-tower}
Let $G$ be a commutative monoid and $H$ an abelian group. For every $n\ge 1$, the powers of $\mathsf{D}_{\langle 2\rangle} = \mathsf{L}+\mathsf{J}_0$ are the higher mixed differences at the origin,
\[
\mathsf{D}_{\langle 2\rangle}^{\,n} = \mathsf{D}_{\langle n+1\rangle}.
\]
\end{proposition}

\begin{proof}
By Lemma~\ref{lem:basic}(ii), every lift of $\mathsf{D}_{\langle 2\rangle}$ splits as $\mathsf{D}_{\langle 2\rangle}^{[i]} = \mathsf{L}^{[i]} + \mathsf{J}_0^{[i]}$, so Lemma~\ref{lem:substitution} shows that every lift of $\mathsf{D}_{\langle 2\rangle}$ carries $\mathsf{D}_{\langle n\rangle}$ to $\mathsf{D}_{\langle n+1\rangle}$. The statement follows by induction along any tower, starting from Proposition~\ref{prop:D2} and using Theorem~\ref{thm:towers}.
\end{proof}

\begin{remark}\label{rem:two-towers}
Let $G$ be a commutative monoid and $H$ an abelian group. The operators $\mathsf{L}$ and $\mathsf{L}+\mathsf{J}_0$ differ by the evaluation at the origin. The powers of the first are the Cauchy differences, and the powers of the second are the mixed differences at the origin. The two families differ only in the term $S=\emptyset$ of (\ref{eq:Dn}), which the Cauchy differences omit, so that $\mathsf{L}^{\,n}[f] = \mathsf{D}_{\langle n+1\rangle}[f] + (-1)^{n}f(0)$, as the comparison of (\ref{eq:Dn}) with the formula of Example~\ref{ex:powers} shows. In particular $\mathsf{F}[f]=\mathsf{D}_{\langle 3\rangle}[f]+f(0)$, so the equation $\mathsf{F}[f]=0$ is Fréchet's equation $\mathsf{D}_{\langle 3\rangle}[f]=0$ together with the condition $f(0)=0$. Over semigroups the relation between the solutions of $\mathsf{F}[f]=0$ and Fréchet's equation is studied in \cite{Stetkaer2017}. The operators $\mathsf{L}$ and $\mathsf{L}+\mathsf{J}_0$ both belong to the pencil of Proposition~\ref{prop:D2}, all of whose members are associative.
\end{remark}

\subsection*{The tower spectrum}
By Theorem~\ref{thm:towers} an associative operator has a single tower of each length. In general the towers of a given length need not coincide, and their number can be counted. For binary operations such a count is known. The associative spectrum of a binary operation is the sequence counting, for each $n$, the distinct term operations induced by the bracketings of $x_1,\dots,x_n$. It measures the failure of the generalised associative law and is bounded by the Catalan numbers \cite{CsakanyWaldhauser,LiebscherWaldhauser2009}. The following notion extends it from operations to operators.

\begin{definition}\label{def:tower-spectrum}
Let $\mathsf{T}\colon \mathcal{H}_1\to \mathcal{H}_2$ be an operator, and for $n\ge 1$ let $t_n(\mathsf{T})$ be the number of distinct towers of length~$n$ over~$\mathsf{T}$. The sequence $(t_n(\mathsf{T}))_{n\ge 1}$ is called the \emph{tower spectrum} of~$\mathsf{T}$.
\end{definition}

For a composition operator $\mu^{*}$ the towers of length $n$ are, by Lemma~\ref{lem:basic}(i), the operators $\nu^{*}$, where $\nu\colon G^{n+1}\to G$ runs over the term operations induced by the bracketings of $x_1,\dots,x_{n+1}$. Every bracketing arises in this way, an innermost product of two adjacent variables being formed by the last lift. For a target with at least two elements two such operators coincide only when the term operations do, so that $t_n(\mu^{*})$ is then the term of index $n+1$ of the associative spectrum of~$\mu$.

\begin{proposition}\label{prop:tower-spectrum}
\leavevmode
\begin{enumerate}
\item[\rm (i)] For every operator $\mathsf{T}\colon \mathcal{H}_1\to \mathcal{H}_2$ and every $n\ge 1$ one has $t_n(\mathsf{T})\le n!$. Moreover, $t_2(\mathsf{T})=1$ if and only if $\mathsf{T}$ is associative, and then $t_n(\mathsf{T})=1$ for every~$n$.
\item[\rm (ii)] If $G$ is a commutative monoid with an element of infinite order and $H$ is a non-trivial abelian group, then $t_n(\Delta)=2^{\,n-1}$ for every~$n$.
\item[\rm (iii)] If $G$ is infinite and $H$ is a free magma on infinitely many generators, with product~$\beta$, then the operator $\mathsf{T}[u](x,y)=\beta(u(x),u(y))$ of Example~\ref{ex:no-exchange} has $t_n(\mathsf{T})=n!$ for every~$n$.
\end{enumerate}
\end{proposition}

\begin{proof}
(i) A tower of length $n$ is given by integers $j_1,\dots,j_n$ with $1\le j_k\le k$ for every~$k$, and there are $n!$ such sequences. The towers of length $2$ are $\mathsf{T}^{[1]}\mathsf{T}$ and $\mathsf{T}^{[2]}\mathsf{T}$, so the equality $t_2(\mathsf{T})=1$ is Definition~\ref{def:associative}, and the last assertion is Theorem~\ref{thm:towers}.

(ii) The multiples of an element of infinite order are pairwise distinct, so the last assertion of Proposition~\ref{prop:delta-towers} applies for every~$n$.

(iii) For $n=1$ there is nothing to prove. Let $n\ge 2$, and let $u$ take pairwise distinct generators as values at pairwise distinct points $x_1,\dots,x_{n+1}$. By (\ref{eq:lift}), the value at $(x_1,\dots,x_{n+1})$ of a tower $\mathsf{T}^{[j_n]}\cdots\mathsf{T}^{[j_1]}[u]$ is the product of two values of the shorter tower $\mathsf{T}^{[j_{n-1}]}\cdots\mathsf{T}^{[j_1]}[u]$, the left one at the point obtained by deleting $x_{j_n+1}$ and the right one at the point obtained by deleting $x_{j_n}$. By induction on~$n$, each of the letters $u(x_1),\dots,u(x_{n+1})$ occurs in this value, and its left factor contains all of them except $u(x_{j_n+1})$. In a free magma a product determines its two factors, so the value determines $j_n$ and then, by a second induction, the whole sequence $j_1,\dots,j_n$. Hence the $n!$ sequences give $n!$ distinct towers.
\end{proof}

The hypothesis on $G$ in Proposition~\ref{prop:tower-spectrum}(ii) cannot be dropped: by Remark~\ref{rem:delta-sharp}, over $G=H=\mathbb{Z}_2$ one has $t_n(\Delta)=2$ for every $n\ge 2$. Question~\ref{q:tower-spectrum} asks which sequences occur as tower spectra.

\section{A classification of the associative integer combinations}\label{sec:sums}

In this section $G$ and $H$ are abelian groups, so that the operators $\mathcal{H}_1\to \mathcal{H}_2$ form an abelian group under pointwise addition. The five basic operators of Section~\ref{sec:lifts}, namely the two translations $\mathsf{E}_{+}$, $\mathsf{E}_{-}$, the embeddings $\mathsf{J}_1$, $\mathsf{J}_2$ and the evaluation at the origin $\mathsf{J}_0$, are additive, because each is a composition operator. They span the integer family
\begin{equation}\label{eq:global-family}
\mathsf{T}_{a,b} := a_{+}\mathsf{E}_{+}+a_{-}\mathsf{E}_{-}+b_1\mathsf{J}_1+b_2\mathsf{J}_2+b_0\mathsf{J}_0,
\end{equation}
that is, $\mathsf{T}_{a,b}[f](x,y) = a_{+} f(x+y)+a_{-} f(x-y)+b_1 f(x)+b_2 f(y)+b_0 f(0)$, with integer coefficients $a=(a_{+},a_{-})$ and $b=(b_0,b_1,b_2)$. Four of the five are the composition operators of associative operations, namely the addition, the two projections and the constant operation with value~$0$, and are associative for that reason. The fifth, $\mathsf{E}_{-}$, is the composition operator of the subtraction, which is associative only when $2G=0$ (Section~\ref{sec:lifts}). The family contains $\mathsf{L}$, $\mathsf{P}$, $\mathsf{Q}$, $\Delta$ and $\mathsf{D}_{\langle 2\rangle}=\mathsf{L}+\mathsf{J}_0$, and Theorem~\ref{thm:global} determines its associative members.

We first treat sums of additive operators. When both summands are associative, the sum is associative exactly when the pair satisfies a bilinear condition, which Definition~\ref{def:compatible} calls compatibility. Example~\ref{ex:incompatible} reads the non-associativity of $\Delta$ as a failure of this condition, and the classification provides two compatible pairs (Corollary~\ref{cor:pencils}). The classification itself is proved by a direct expansion, which also covers the members involving $\mathsf{E}_{-}$.

\begin{proposition}\label{prop:sum}
Let $\mathsf{T},\mathsf{S}\colon \mathcal{H}_1\to \mathcal{H}_2$ be additive and associative. Then $\mathsf{T}+\mathsf{S}$ is associative if and only if
\begin{equation}\label{eq:mixed}
\mathsf{T}^{[1]}\mathsf{S} + \mathsf{S}^{[1]}\mathsf{T} \;=\; \mathsf{T}^{[2]}\mathsf{S} + \mathsf{S}^{[2]}\mathsf{T}.
\end{equation}
\end{proposition}

\begin{proof}
By Lemma~\ref{lem:basic}(ii), $(\mathsf{T}+\mathsf{S})^{[i]} = \mathsf{T}^{[i]}+\mathsf{S}^{[i]}$, and by Lemma~\ref{lem:basic}(iii) the four lifts involved are additive. Hence, for $i=1,2$,
\[
(\mathsf{T}+\mathsf{S})^{[i]}(\mathsf{T}+\mathsf{S})
= \mathsf{T}^{[i]}\mathsf{T} + \mathsf{S}^{[i]}\mathsf{S} + \mathsf{T}^{[i]}\mathsf{S} + \mathsf{S}^{[i]}\mathsf{T},
\]
and when the two instances are subtracted the terms $\mathsf{T}^{[i]}\mathsf{T}$ and $\mathsf{S}^{[i]}\mathsf{S}$ cancel, by the associativity of $\mathsf{T}$ and of~$\mathsf{S}$. What remains is the difference of the two sides of~(\ref{eq:mixed}).
\end{proof}

\begin{definition}\label{def:compatible}
Two additive associative operators $\mathsf{T},\mathsf{S}\colon \mathcal{H}_1\to \mathcal{H}_2$ are called \emph{compatible} if they satisfy (\ref{eq:mixed}).
\end{definition}

\begin{remark}\label{rem:compatible-combinations}
Two additive associative operators $\mathsf{T}$ and $\mathsf{S}$ are compatible if and only if every integer combination $m\mathsf{T}+n\mathsf{S}$ is associative. Indeed, each multiple of an additive associative operator is again additive and associative, by Lemma~\ref{lem:basic}(ii) and~(iii), so Proposition~\ref{prop:sum} applies to the pair $(m\mathsf{T},n\mathsf{S})$.
\end{remark}

\begin{remark}\label{rem:compatible}
Condition (\ref{eq:mixed}) is the operator analogue of the \emph{compatibility} of two associative products, the requirement that their sum, equivalently every linear combination, be associative. For bilinear products this notion, governed by the vanishing of the mixed associator, is studied by Odesskii and Sokolov \cite{OdesskiiSokolov}. For operators that are not additive, the expansion in the proof of Proposition~\ref{prop:sum} is not available, and the associativity of a sum has to be examined case by case.
\end{remark}

\begin{example}\label{ex:incompatible}
When $G$ and $H$ are non-trivial, the pair $(\mathsf{E}_{+},-\mathsf{J}_2)$ is not compatible. Indeed, the four mixed terms are computed directly from (\ref{eq:lift}):
\[
\begin{gathered}
\mathsf{E}_{+}^{[1]}\mathsf{J}_2[f](x,y,z) = f(z), \qquad \mathsf{J}_2^{[1]}\mathsf{E}_{+}[f](x,y,z) = f(y+z),\\
\mathsf{E}_{+}^{[2]}\mathsf{J}_2[f](x,y,z) = f(y+z), \qquad \mathsf{J}_2^{[2]}\mathsf{E}_{+}[f](x,y,z) = f(x+z),
\end{gathered}
\]
so the two sides of (\ref{eq:mixed}) for $(\mathsf{E}_{+},-\mathsf{J}_2)$ differ by exactly $f(x+z)-f(z)$, which is the failure of associativity of $\Delta=\mathsf{E}_{+}-\mathsf{J}_2$ recorded in (\ref{eq:delta-failure}).
\end{example}

The Cauchy difference $\mathsf{L}=\mathsf{E}_{+}-\mathsf{J}_1-\mathsf{J}_2$ shows that a combination of three associative operators can be associative although, for non-trivial $G$ and $H$, no two of them are compatible. Indeed, the two sides of (\ref{eq:mixed}) differ by $f(x+z)-f(z)$ for the pair $(\mathsf{E}_{+},-\mathsf{J}_2)$, by $f(x)-f(x+z)$ for $(\mathsf{E}_{+},-\mathsf{J}_1)$ and by $f(z)-f(x)$ for $(-\mathsf{J}_1,-\mathsf{J}_2)$, and the three differences cancel. The proof of the classification uses the following independence lemma.

\begin{lemma}\label{lem:independence}
Let $r,s\ge 1$ and let $p\in G^{s}$. Let $\mu_1,\dots,\mu_r\colon G^{s}\to G$ be maps and let the integers $m_1,\dots,m_r$ satisfy $\sum_k m_k f(\mu_k(p)) = 0$ for all $f\in \mathcal{H}_1$. Then, for every $v\in G$,
\[
\Bigl(\,\sum_{k\,:\,\mu_k(p)=v} m_k\Bigr) H = 0.
\]
In particular, if the values $\mu_1(p),\dots,\mu_r(p)$ are pairwise distinct, then $m_kH=0$ for every~$k$, and if moreover $H$ has infinite exponent then $m_1=\dots=m_r=0$.
\end{lemma}

\begin{proof}
Fix $v\in G$ and $c\in H$, and let $f$ take the value $c$ at $v$ and $0$ elsewhere. The hypothesis then reads $\bigl(\sum_{k:\mu_k(p)=v}m_k\bigr)c=0$, and $c\in H$ was arbitrary.
\end{proof}

The classification is stated in terms of the seven integers
\begin{equation}\label{eq:seven}
\begin{gathered}
N_1 := a_{-}(a_{+}{-}a_{-}), \qquad N_2 := a_{+}(b_1{-}b_2), \qquad N_3 := a_{-}(b_1{-}b_2),\\
N_4 := b_1(a_{+}{+}a_{-}{+}b_2)+b_0(a_{+}{+}a_{-}{+}b_1), \qquad
N_5 := b_2(a_{+}{+}a_{-}{+}b_1)+b_0(a_{+}{+}b_2),\\
N_6 := b_0\,a_{-}, \qquad N_7 := b_0(b_1{-}b_2),
\end{gathered}
\end{equation}
attached to the coefficients $a=(a_{+},a_{-})$ and $b=(b_0,b_1,b_2)$. They are the coefficients that remain in the expansion of $\mathsf{T}_{a,b}^{[1]}\mathsf{T}_{a,b}-\mathsf{T}_{a,b}^{[2]}\mathsf{T}_{a,b}$ carried out in the proof of Theorem~\ref{thm:global}.

\begin{theorem}[classification]\label{thm:global}
Let $G$ and $H$ be abelian groups. Suppose $G$ has an element $g$ with $2g\neq 0$, and suppose that $H$ has infinite exponent, that is, that no non-zero integer annihilates it. Let $a_{+},a_{-},b_0,b_1,b_2\in\mathbb{Z}$. The following are equivalent:
\begin{enumerate}
\item[\rm (i)] $\mathsf{T}_{a,b}$ is associative;
\item[\rm (ii)] the seven integers (\ref{eq:seven}) all vanish;
\item[\rm (iii)] $\mathsf{T}_{a,b}$ is an integer multiple of one of $\mathsf{J}_1$, $\mathsf{J}_2$, $\mathsf{P}$, $\mathsf{Q}$, or lies in one of the two pencils $\{m\mathsf{L}+n\mathsf{J}_0\}_{m,n\in\mathbb{Z}}$ and $\{m\mathsf{E}_{+}+n(\mathsf{J}_0-\mathsf{J}_1-\mathsf{J}_2)\}_{m,n\in\mathbb{Z}}$.
\end{enumerate}
The implication {\rm (ii)}$\Rightarrow${\rm (i)} holds for all abelian groups $G$ and $H$.
\end{theorem}

\begin{proof}
Write $\mathsf{T}:=\mathsf{T}_{a,b}$ and $\psi:=\mathsf{T}[f]$. Expanding the two iterates along (\ref{eq:lift}) and cancelling the terms that agree, the difference $\mathsf{T}^{[1]}\mathsf{T}[f]-\mathsf{T}^{[2]}\mathsf{T}[f]$ equals
\begin{multline}\label{eq:global-difference}
N_1\bigl[f(x-y+z)-f(x-y-z)\bigr] + N_2\,f(x+z) + N_3\,f(x-z)\\
- N_4\,f(x) + N_5\,f(z) + N_6\,f(-z) + N_7\,f(0).
\end{multline}
Indeed, the two iterates are
\[
\begin{gathered}
a_{+}\psi(x{+}y,z)+a_{-}\psi(x{-}y,z)+b_1\psi(x,z)+b_2\psi(y,z)+b_0\psi(0,z),\\
a_{+}\psi(x,y{+}z)+a_{-}\psi(x,y{-}z)+b_1\psi(x,y)+b_2\psi(x,z)+b_0\psi(x,0),
\end{gathered}
\]
and expanding each of the ten evaluations of $\psi$ through (\ref{eq:global-family}) produces integer combinations of the fifteen maps
\[
x\pm y\pm z,\quad x\pm y,\quad x\pm z,\quad y\pm z,\quad x,\quad y,\quad z,\quad -z,\quad 0
\]
composed with~$f$. Collecting coefficients leaves exactly (\ref{eq:global-difference}). The collection is a formal identity between integer combinations of maps $G^3\to G$, valid whatever coincidences occur among the values of the maps that cancel. If the seven integers (\ref{eq:seven}) vanish, then (\ref{eq:global-difference}) vanishes identically and $\mathsf{T}$ is associative. This proves (ii)$\Rightarrow$(i) for all abelian groups $G$ and $H$, which is the last assertion of the theorem.

Assume now (i), and fix $g\in G$ with $2g\neq 0$, so that $g$ differs from $0$, from $-g$ and from $2g$. Only the first bracket of (\ref{eq:global-difference}) involves~$y$, and for fixed $x$ the substitution $y\mapsto x-s$ makes $s$ run over the whole of~$G$, turning that bracket into $N_1[f(s+z)-f(s-z)]$. The vanishing of (\ref{eq:global-difference}) therefore forces the latter to be independent of~$s$. Taking $z=g$ and comparing $s=0$ with $s=g$ gives, for every $f\in \mathcal{H}_1$,
\[
N_1\bigl[f(g)-f(-g)-f(2g)+f(0)\bigr] = 0,
\]
and Lemma~\ref{lem:independence} at the value $g$, together with the hypothesis on~$H$, gives $N_1=0$.

The six remaining terms of (\ref{eq:global-difference}) do not involve~$y$. At $z=0$ and $x=g$ they reduce to
\[
(N_2+N_3-N_4)\,f(g) + (N_5+N_6+N_7)\,f(0),
\]
and at $x=0$ and $z=g$ to
\[
(N_2+N_5)\,f(g) + (N_3+N_6)\,f(-g) + (N_7-N_4)\,f(0).
\]
The values involved are pairwise distinct in each case, so Lemma~\ref{lem:independence} yields $N_5=-N_2$, $N_6=-N_3$ and $N_4=N_7=N_2+N_3$. At $x=z=g$ the six terms become
\[
N_2\,f(2g) - (2N_2+N_3)\,f(g) - N_3\,f(-g) + (N_2+2N_3)\,f(0).
\]
If $3g\neq 0$, the four values $2g$, $g$, $-g$, $0$ are pairwise distinct and Lemma~\ref{lem:independence} gives $N_2=N_3=0$ at once. If $3g=0$, then $2g=-g$, and the three values $-g$, $g$, $0$ carry the coefficients $N_2-N_3$, $-(2N_2+N_3)$ and $N_2+2N_3$. The first two give $N_2=N_3$ and $3N_2=0$, hence again $N_2=N_3=0$. All seven integers (\ref{eq:seven}) therefore vanish, which proves (i)$\Rightarrow$(ii).

The coefficients of each family in (iii) satisfy (ii), as a direct check shows. Every operator in (iii) is therefore associative, by the implication (ii)$\Rightarrow$(i), which gives (iii)$\Rightarrow$(i). For (ii)$\Rightarrow$(iii) we use that $\mathbb{Z}$ is an integral domain. There are three cases.

\emph{Case $a_{-}\neq 0$.} From $N_1=0$ we get $a_{+}=a_{-}$, from $N_6=0$ we get $b_0=0$, and from $N_3=0$ we get $b_1=b_2=:d$. Then $N_4$ becomes $d(2a_{+}+d)$, so either $d=0$, giving the multiples of $\mathsf{P}=\mathsf{E}_{+}+\mathsf{E}_{-}$, or $d=-2a_{+}$, giving the multiples of $\mathsf{Q}$.

\emph{Case $a_{-}=0$ and $b_1\neq b_2$.} Here $N_2=0$ and $N_7=0$ force $a_{+}=b_0=0$, and the two conditions $N_4=N_5=0$ reduce to $b_1b_2=0$, giving the multiples of $\mathsf{J}_1$ and of $\mathsf{J}_2$.

\emph{Case $a_{-}=0$ and $b_1=b_2=:d$.} The two conditions $N_4=0$ and $N_5=0$ both read
\[
d^2+d(a_{+}+b_0)+b_0\,a_{+} = (d+a_{+})(d+b_0) = 0,
\]
so $d=-a_{+}$, giving the pencil $m\mathsf{L}+n\mathsf{J}_0$ with $m=a_{+}$, $n=b_0$, or $d=-b_0$, giving the pencil $m\mathsf{E}_{+}+n(\mathsf{J}_0-\mathsf{J}_1-\mathsf{J}_2)$ with $m=a_{+}$, $n=b_0$.
\end{proof}

The classification singles out the operators of Sections \ref{sec:lifts} and~\ref{sec:towers} inside the family (\ref{eq:global-family}).

\begin{corollary}\label{cor:detect}
Under the hypotheses of Theorem~\ref{thm:global}, the Cauchy difference is, up to an integer multiple, the unique associative member of the subfamily $a_{+}\mathsf{E}_{+}+b_1\mathsf{J}_1+b_2\mathsf{J}_2$ with $a_{+},b_1,b_2$ all non-zero. The Jordan--von Neumann operator is, up to a multiple, the unique associative member of $a_{+}\mathsf{E}_{+}+a_{-}\mathsf{E}_{-}+b_1\mathsf{J}_1+b_2\mathsf{J}_2$ with $a_{+},a_{-},b_1,b_2$ all non-zero, and the d'Alembert operator is, up to a multiple, the unique associative member of $a_{+}\mathsf{E}_{+}+a_{-}\mathsf{E}_{-}+b_0\mathsf{J}_0$ with $a_{+},a_{-}$ both non-zero. The intersection of the two pencils consists of the integer multiples of the mixed difference at the origin $\mathsf{D}_{\langle 2\rangle}$. Finally, the finite difference $\Delta=\mathsf{E}_{+}-\mathsf{J}_2$ is not associative.
\end{corollary}

\begin{proof}
Under the hypotheses of Theorem~\ref{thm:global} the five spanning operators are independent over $\mathbb{Z}$, so that a member of (\ref{eq:global-family}) determines its coefficients. Indeed, if $\mathsf{T}_{a,b}=0$, then evaluating at $(0,g)$, at $(g,g)$ and at $(g,0)$, with $2g\neq 0$, and applying Lemma~\ref{lem:independence} gives in turn $a_{-}=0$, then $a_{+}=b_0=0$, then $b_1=b_2=0$. All the statements are now read off from Theorem~\ref{thm:global}(iii). For the intersection of the pencils, $m\mathsf{L}+n\mathsf{J}_0 = m'\mathsf{E}_{+}+n'(\mathsf{J}_0-\mathsf{J}_1-\mathsf{J}_2)$ forces, on comparing the coefficients of $\mathsf{E}_{+}$, $\mathsf{J}_1$ and $\mathsf{J}_0$, first $m=m'$, then $m=n'$, then $n=n'$. Hence $m=n$, and the operator is $m(\mathsf{L}+\mathsf{J}_0)=m\mathsf{D}_{\langle 2\rangle}$.
\end{proof}

\begin{corollary}\label{cor:pencils}
Let $G$ and $H$ be abelian groups. The pairs $(\mathsf{L},\mathsf{J}_0)$ and $(\mathsf{E}_{+},\,\mathsf{J}_0-\mathsf{J}_1-\mathsf{J}_2)$ are compatible.
\end{corollary}

\begin{proof}
The coefficients of every integer combination of either pair satisfy condition (ii) of Theorem~\ref{thm:global}. By the last assertion of Theorem~\ref{thm:global}, every such combination is associative, and so are in particular the four operators themselves. By Remark~\ref{rem:compatible-combinations} the two pairs are therefore compatible.
\end{proof}

\begin{remark}[the two hypotheses]\label{rem:sharp}
Neither of the two hypotheses of Theorem~\ref{thm:global} can be dropped. If $2G=0$, the two translations coincide, $\mathsf{E}_{-}=\mathsf{E}_{+}$, and the operator $\mathsf{E}_{+}-\mathsf{E}_{-}=0$ is associative although $N_1=-2$. If $mH=0$ for some $m\neq 0$, the operator $m\mathsf{J}_1+m\mathsf{J}_2=0$ is associative although $N_4=N_5=m^2$. For $m=2$ the operator $\mathsf{E}_{+}+\mathsf{J}_1+\mathsf{J}_2$, which coincides with $\mathsf{L}$ when $2H=0$, is associative although $N_4=N_5=2$. In each case it is the implication (i)$\Rightarrow$(ii) that fails.

The hypothesis on $G$ is needed for the criterion~(ii) only: if $H$ has infinite exponent, then (i) and (iii) are equivalent for every abelian group~$G$. Indeed, (iii) implies (i) by the last assertion of Theorem~\ref{thm:global}, since the coefficients of each family in (iii) satisfy~(ii). By the theorem itself it remains to derive (iii) from (i) when $2G=0$. Then $\mathsf{E}_{-}=\mathsf{E}_{+}$, so that $\mathsf{T}_{a,b}=\mathsf{T}_{a',b}$ with $a':=(a_{+}+a_{-},0)$. If $G=0$, every member of the family (\ref{eq:global-family}) is an integer multiple of~$\mathsf{J}_1$. Otherwise let $g\neq 0$, let $\mathsf{T}_{a,b}$ be associative, and let $N_1,\dots,N_7$ now denote the integers (\ref{eq:seven}) attached to $(a',b)$. Then $N_1=N_3=N_6=0$, so that (\ref{eq:global-difference}), written for $(a',b)$, reduces to $N_2\,f(x+z)-N_4\,f(x)+N_5\,f(z)+N_7\,f(0)$, and this vanishes for all $f$, $x$ and~$z$. At $(x,z)=(g,0)$, at $(0,g)$ and at $(g,g)$, where $g+g=0$, Lemma~\ref{lem:independence} gives in turn $N_2=N_4$ and $N_5=-N_7$, then $N_5=-N_2$ and $N_7=N_4$, then $N_7=-N_2$. Hence $2N_2=0$, and the seven integers attached to $(a',b)$ vanish. The proof of the implication (ii)$\Rightarrow$(iii) uses no hypothesis on $G$ or on~$H$, so $\mathsf{T}_{a',b}$ is one of the operators in~(iii).
\end{remark}

\section{Symmetry}\label{sec:symmetry}

Many of the operators of Section~\ref{sec:lifts} take symmetric functions as values. This section relates symmetry to associativity. Symmetry also distinguishes two of the image problems of Section~\ref{sec:image}: the image of the Cauchy difference over a commutative semigroup consists of symmetric functions, whereas the image of the Jordan--von Neumann operator does not when $2G\neq 0$ and $H$ is non-trivial.

\begin{definition}\label{def:symmetric}
A function $\psi\in \mathcal{H}_n$ is called \emph{symmetric} if it is invariant under every permutation of its $n$ arguments. An operator $\mathsf{T}\colon \mathcal{H}_1\to \mathcal{H}_2$ is called \emph{symmetric} if $\mathsf{T}[f]$ is symmetric for every $f\in \mathcal{H}_1$, that is, if $\mathsf{T}[f](x,y) = \mathsf{T}[f](y,x)$ for all $f\in \mathcal{H}_1$ and $x,y\in G$.
\end{definition}

\begin{proposition}\label{prop:symmetric-model}
Let $\mu\colon G^2\to G$ and let $H$ have at least two elements. The composition operator $\mu^{*}$ is symmetric if and only if the operation $\mu$ is commutative.
\end{proposition}

\begin{proof}
If $\mu$ is commutative, then $\mu^{*}[f](x,y)=f(\mu(x,y))=f(\mu(y,x))=\mu^{*}[f](y,x)$. Conversely, if $\mu(x,y)\neq\mu(y,x)$ for some $x,y$, choose $c_0\neq c_1$ in $H$ and let $f$ take the value $c_0$ at $\mu(x,y)$ and $c_1$ elsewhere. Then $\mu^{*}[f]$ takes different values at $(x,y)$ and $(y,x)$.
\end{proof}

\begin{example}\label{ex:symmetric}
Over a commutative semigroup, $\mathsf{E}$ and $\mathsf{L}$ are symmetric, and over a commutative monoid so are $\mathsf{J}_0$ and $\mathsf{D}_{\langle 2\rangle}=\mathsf{L}+\mathsf{J}_0$. Indeed, the values $f(x+y)$, $f({x+y})-f(x)-f(y)$, $f(0)$ and $f(x+y)-f(x)-f(y)+f(0)$ are unchanged when $x$ and $y$ are interchanged. By contrast, let $G$ have at least two elements and let $H$ be a non-trivial abelian group. The embeddings $\mathsf{J}_1$ and $\mathsf{J}_2$ are then not symmetric. Over a commutative semigroup the difference $\Delta$ is not symmetric either, and neither are $\mathsf{P}$ and $\mathsf{Q}$ when $2G\neq 0$, that is, when $x-y\neq y-x$ for some $x,y$. The term $f(x-y)$ is the one that is not symmetric. Symmetry and associativity are independent. The operator $\mathsf{J}_1$ is associative and not symmetric, while the composition operator of a commutative non-associative operation is symmetric and not associative, and so is, under the hypotheses of Theorem~\ref{thm:global}, the operator $\mathsf{E}_{+}+\mathsf{J}_0$, whose coefficients violate condition~(ii).
\end{example}

\begin{theorem}[symmetric towers]\label{thm:symmetric-towers}
Let $\mathsf{T}$ be associative and symmetric. Then, for every $n\ge 1$, the power $\mathsf{T}^{\,n}[f]$ is a symmetric function of its $n+1$ arguments, for every $f\in \mathcal{H}_1$.
\end{theorem}

\begin{proof}
For $n=1$ the assertion is the symmetry of $\mathsf{T}$, so let $n\ge 2$ and fix $i$ with $1\le i\le n$. By Theorem~\ref{thm:towers}, $\mathsf{T}^{\,n} = \mathsf{T}^{[i]}\,\mathsf{T}^{\,n-1}$, so, by (\ref{eq:lift}),
\[
\mathsf{T}^{\,n}[f](x_1,\dots,x_{n+1}) = \mathsf{T}\bigl[\mathsf{T}^{\,n-1}[f](x_1,\dots,x_{i-1},\,\cdot\,,x_{i+2},\dots,x_{n+1})\bigr](x_i,x_{i+1}).
\]
The partial map in the bracket does not involve $x_i$ or $x_{i+1}$, and $\mathsf{T}[u]$ is a symmetric function for every $u\in \mathcal{H}_1$. Hence the displayed value is invariant under the transposition of $x_i$ and $x_{i+1}$. As $i$ was arbitrary, $\mathsf{T}^{\,n}[f]$ is invariant under all the adjacent transpositions $(x_i,x_{i+1})$, $1\le i\le n$, and these generate the full symmetric group on the $n+1$ arguments.
\end{proof}

\begin{remark}\label{rem:symmetric-towers}
The proof uses Theorem~\ref{thm:towers} for every slot. The invariance under the transposition of $x_i$ and $x_{i+1}$ comes from the representation $\mathsf{T}^{[i]}\mathsf{T}^{\,n-1}$ of $\mathsf{T}^{\,n}$, and these representations are available for all $i$ because all towers coincide. Theorem~\ref{thm:symmetric-towers} gives, for instance, a second proof of the symmetry of the mixed differences at the origin. By Proposition~\ref{prop:origin-tower}, $\mathsf{D}_{\langle n+1\rangle}$ is the $n$-th power of the symmetric associative operator $\mathsf{D}_{\langle 2\rangle}$, so its values are symmetric in all $n+1$ increments. In Section~\ref{sec:towers} the same fact follows from the commutation of the family $(\Delta_t)$.
\end{remark}

\begin{proposition}\label{prop:symmetric-classification}
Under the hypotheses of Theorem~\ref{thm:global}, an operator $\mathsf{T}_{a,b}$ of the family (\ref{eq:global-family}) is symmetric and associative if and only if it lies in one of the two pencils of Theorem~\ref{thm:global}{\rm (iii)}.
\end{proposition}

\begin{proof}
Every member of either pencil is associative by Theorem~\ref{thm:global}, and symmetric because its formula $m f(x+y)+d f(x)+d f(y)+n f(0)$, with $d=-m$ or $d=-n$, is invariant under the transposition of $x$ and~$y$. Conversely, by Theorem~\ref{thm:global} it suffices to rule out the non-zero multiples of $\mathsf{J}_1$, $\mathsf{J}_2$, $\mathsf{P}$ and $\mathsf{Q}$. Let $g\in G$ satisfy $2g\neq 0$, so that $g\neq 0$ and $-g\neq g$, and, for $c\in H$, let $f_c$ take the value $c$ at $g$ and $0$ elsewhere. For $m\mathsf{J}_1$, symmetry at $(g,0)$ gives $mf_c(g)=mf_c(0)$, that is, $mc=0$ for every $c\in H$, so $m=0$ by the hypothesis on $H$. The same argument applies to $m\mathsf{J}_2$. For $m\mathsf{P}$, symmetry at $(g,0)$ reads $2mf_c(g)=m\bigl[f_c(g)+f_c(-g)\bigr]$, that is, $mf_c(g)=mf_c(-g)$, and the choice of $f_c$ gives $mc=0$ for every $c\in H$, so $m=0$. For $m\mathsf{Q}$ the same comparison at $(g,0)$, after cancelling the common terms, gives again $mf_c(g)=mf_c(-g)$ and $m=0$.
\end{proof}

\begin{remark}[quasi-symmetry]\label{rem:quasi-symmetry}
If $\mathsf{T}$ is symmetric, every function in $\operatorname{Im}\mathsf{T}$ is symmetric, and if $\mathsf{T}$ is moreover associative, every function in $\operatorname{Im}\mathsf{T}^{\,n}$ is symmetric, by Theorem~\ref{thm:symmetric-towers}. This is the source of the symmetry requirement in the classical description of $\operatorname{Im}\mathsf{L}$ recalled in Section~\ref{sec:image}. Let now $G$ and $H$ be abelian groups. If $2G\neq 0$ and $H$ is non-trivial, the Jordan--von Neumann operator is not symmetric. Every $\psi$ in its image satisfies instead the \emph{quasi-symmetry} identity $\psi(x,y)-\psi(y,x)=\psi(0,0)-\psi(0,x-y)$, as a direct expansion of both sides shows. Suppose that $H$ is uniquely $2$-divisible. Every $f$ then splits as $f=f_{\mathrm{e}}+f_{\mathrm{o}}$ with $f_{\mathrm{e}}$ even and $f_{\mathrm{o}}$ odd, and $\mathsf{Q}[f]=\mathsf{Q}[f_{\mathrm{e}}]+\mathsf{Q}[f_{\mathrm{o}}]$. The function $\mathsf{Q}[f_{\mathrm{e}}]$ is symmetric, so the lack of symmetry comes from the odd part of $f$ alone:
\[
\mathsf{Q}[f_{\mathrm{o}}](x,y)-\mathsf{Q}[f_{\mathrm{o}}](y,x) = f_{\mathrm{o}}(x-y)-f_{\mathrm{o}}(y-x) = 2f_{\mathrm{o}}(x-y).
\]
\end{remark}

\section{Cocycle identities and the image problem}\label{sec:image}

We return to the inverse problem of Section~\ref{sec:intro}: for which $\psi\in \mathcal{H}_2$ is the equation $\mathsf{T}[f]=\psi$ solvable? For the Cauchy difference the classical answer describes $\operatorname{Im}\mathsf{L}$ as the set of symmetric solutions of the \emph{cocycle equation} $\psi(x,y)+\psi(x+y,z)=\psi(y,z)+\psi(x,y+z)$. This goes back to Kurepa, who solved the equation for differentiable real functions \cite{Kurepa1956}, and to Erdős, who identified the symmetric solutions as the coboundaries $f(x+y)-f(x)-f(y)$, for real functions and by way of the divisibility of the additive group of the reals \cite{Erdos1959}. Under hypotheses of that kind the general solution of the cocycle equation is such a coboundary plus a skew-symmetric biadditive map. Some hypothesis is needed. Over $G=H=\mathbb{Z}_3$ the factor set of the non-split extension $\mathbb{Z}_9$ is a symmetric cocycle that is not a coboundary, while the only skew-symmetric biadditive map $\mathbb{Z}_3\times\mathbb{Z}_3\to\mathbb{Z}_3$ is zero.

For real functions, Hosszú obtained the corresponding splitting of an arbitrary solution \cite{Hosszu1963}, and the statement was then established by Aczél on abelian groups with uniquely divisible targets \cite{Aczel1965} and by Jessen, Karpf and Thorup for torsion-free divisible modules \cite{JessenKarpfThorup1968}. Hosszú stated it for cancellative abelian semigroups with uniquely divisible targets \cite{Hosszu1971}. Davison and Ebanks reported a gap in his proof and gave a new one, which does not use the uniqueness of division. In its final form, that every symmetric cocycle on a cancellative abelian monoid with values in a divisible abelian group is a coboundary, the statement is due to them \cite{DavisonEbanks1995}. See also \cite{EilenbergMacLane,Prunescu2007}. When the domain is a periodic abelian semigroup and the target is uniquely divisible, symmetry is automatic and \emph{every} cocycle is a coboundary \cite{DavisonEbanks1998}. For the Cauchy difference over a periodic abelian group, a finite one for instance, and with such a target, the cocycle identity is therefore the complete answer to the inverse problem.

For the Jordan--von Neumann operator over the reals, a cocycle identity is necessary, and for a datum of class $C^{2}$ it is also sufficient, with the solution given by a closed-form integral \cite{PopaRusJvN}. Over arbitrary abelian groups the image also satisfies the quasi-symmetry identity of Remark~\ref{rem:quasi-symmetry}. This section shows what associativity contributes to these problems for an arbitrary operator, namely the cocycle identities. Symmetry was treated in Section~\ref{sec:symmetry}.

\begin{proposition}[cocycle identity]\label{prop:cocycle-complex}
An operator $\mathsf{T}\colon \mathcal{H}_1\to \mathcal{H}_2$ is associative if and only if every $\psi\in\operatorname{Im}\mathsf{T}$ satisfies
\begin{equation}\label{eq:cocycle-cond}
\mathsf{T}^{[1]}[\psi] \;=\; \mathsf{T}^{[2]}[\psi].
\end{equation}
\end{proposition}

\begin{proof}
The elements of $\operatorname{Im}\mathsf{T}$ are the functions $\psi=\mathsf{T}[f]$ with $f\in \mathcal{H}_1$, and for such a $\psi$ the two sides of (\ref{eq:cocycle-cond}) are $\mathsf{T}^{[1]}\mathsf{T}[f]$ and $\mathsf{T}^{[2]}\mathsf{T}[f]$. They coincide for every $f$ exactly when $\mathsf{T}$ is associative, by Definition~\ref{def:associative}.
\end{proof}

Proposition~\ref{prop:cocycle-complex} is Definition~\ref{def:associative} read on the image. It is in this form that associativity enters the inverse problem: for an associative $\mathsf{T}$, identity (\ref{eq:cocycle-cond}) is a necessary condition on the datum $\psi$ of the equation $\mathsf{T}[f]=\psi$, and it is called the \emph{cocycle identity} of~$\mathsf{T}$. When $H$ is an abelian group, it says that the composite of $\mathsf{T}$ with $\mathsf{T}^{[1]}-\mathsf{T}^{[2]}$ is zero. If $\mathsf{T}$ is moreover additive, the sequence $\mathcal{H}_1\xrightarrow{\ \mathsf{T}\ }\mathcal{H}_2\xrightarrow{\ \mathsf{T}^{[1]}-\mathsf{T}^{[2]}\ }\mathcal{H}_3$ is then a complex of abelian groups. The instances of (\ref{eq:cocycle-cond}) for the operators of this paper are the following:
\begin{enumerate}
\item[\rm (i)] for $\mathsf{T}=\mathsf{L}$, condition (\ref{eq:cocycle-cond}) reads, after cancelling the common term $-\psi(x,z)$ from both sides,
\[
\psi(x,y)+\psi(x+y,z) \;=\; \psi(y,z)+\psi(x,y+z),
\]
the cocycle equation of the Cauchy difference, in the terminology of \cite{DavisonEbanks1995,DavisonEbanks1998};
\item[\rm (ii)] for $\mathsf{T}=\mathsf{Q}$, it reads, after cancelling the common term $-2\psi(x,z)$ from both sides,
\[
\psi(x{+}y,z)+\psi(x{-}y,z)+2\psi(x,y) \;=\; \psi(x,y{+}z)+\psi(x,y{-}z)+2\psi(y,z),
\]
the cocycle identity of \cite{PopaRusJvN};
\item[\rm (iii)] for $\mathsf{T}=\mathsf{P}$, it reads $\psi(x{+}y,z)+\psi(x{-}y,z) = \psi(x,y{+}z)+\psi(x,y{-}z)$, a cocycle identity of d'Alembert type;
\item[\rm (iv)] for an associative $\mu$ and the operator $\mu^{*}$ it reads $\psi(\mu(x,y),z)=\psi(x,\mu(y,z))$; Theorem~\ref{thm:monoid-image} and Proposition~\ref{prop:bands} show that on monoids and on bands this identity characterises the image;
\item[\rm (v)] for a weighted composition operator $\omega\mu^{*}$ with $\mu$ associative and $\omega$ satisfying (\ref{eq:mult-cocycle}), it reads $\omega(x,y)\,\psi(\mu(x,y),z)=\omega(y,z)\,\psi(x,\mu(y,z))$, which is (\ref{eq:weighted-sincov}), the weighted form of the equation in~(iv). For the twisted translation $\omega\mathsf{E}$ it reads
\[
\omega(x,y)\,\psi(x+y,z)=\omega(y,z)\,\psi(x,y+z).
\]
\end{enumerate}

The identity extends to the powers.

\begin{corollary}\label{cor:higher-cocycle}
Let $\mathsf{T}$ be associative and let $n\ge 1$. Then every $\psi\in\operatorname{Im}\mathsf{T}^{\,n}$ satisfies
\[
\mathsf{T}^{[1]}[\psi] = \mathsf{T}^{[2]}[\psi] = \dots = \mathsf{T}^{[n+1]}[\psi].
\]
\end{corollary}

\begin{proof}
If $\psi=\mathsf{T}^{\,n}[f]$, then $\mathsf{T}^{[i]}[\psi] = \mathsf{T}^{[i]}\,\mathsf{T}^{\,n}[f] = \mathsf{T}^{\,n+1}[f]$ for every $i$, by (\ref{eq:power-rec}).
\end{proof}

For a symmetric datum the chain of identities reduces to the first one.

\begin{proposition}[symmetric reduction]\label{prop:symmetric-reduction}
Let $\mathsf{T}\colon \mathcal{H}_1\to \mathcal{H}_2$ be an operator, let $n\ge 1$, and let $\psi\in \mathcal{H}_{n+1}$ be symmetric. If $\mathsf{T}^{[1]}[\psi]=\mathsf{T}^{[2]}[\psi]$, then
\[
\mathsf{T}^{[1]}[\psi] = \mathsf{T}^{[2]}[\psi] = \dots = \mathsf{T}^{[n+1]}[\psi].
\]
\end{proposition}

\begin{proof}
Put
\[
\Lambda(a,b;z_1,\dots,z_n) := \mathsf{T}\bigl[\psi(\,\cdot\,,z_1,\dots,z_n)\bigr](a,b),\qquad a,b,z_1,\dots,z_n\in G.
\]
Since $\psi$ is symmetric, its $i$-th partial map is unchanged when the active variable is moved to the first place, so (\ref{eq:lift}) reads
\begin{equation}\label{eq:sym-red}
\mathsf{T}^{[i]}[\psi](x_1,\dots,x_{n+2}) = \Lambda\bigl(x_i,x_{i+1};\,x_1,\dots,x_{i-1},x_{i+2},\dots,x_{n+2}\bigr)
\end{equation}
for $1\le i\le n+1$, and $\Lambda$ is, for the same reason, a symmetric function of $z_1,\dots,z_n$. Read through (\ref{eq:sym-red}) with $i=1$ and $i=2$, the hypothesis says that
\[
\Lambda(y_1,y_2;y_3,y_4,\dots,y_{n+2}) = \Lambda(y_2,y_3;y_1,y_4,\dots,y_{n+2})
\]
for all $y_1,\dots,y_{n+2}\in G$. Fix $i$ with $1\le i\le n$ and $x_1,\dots,x_{n+2}\in G$, and write $w$ for the tuple $(x_1,\dots,x_{i-1},x_{i+3},\dots,x_{n+2})$ of the variables outside the block $x_i,x_{i+1},x_{i+2}$. Applied to $y_1=x_i$, $y_2=x_{i+1}$, $y_3=x_{i+2}$ and $(y_4,\dots,y_{n+2})=w$, the last identity reads
\[
\Lambda(x_i,x_{i+1};\,x_{i+2},w) = \Lambda(x_{i+1},x_{i+2};\,x_i,w).
\]
Since $\Lambda$ is a symmetric function of the arguments after the semicolon, these can be reordered, and by (\ref{eq:sym-red}) the two sides are then the values of $\mathsf{T}^{[i]}[\psi]$ and $\mathsf{T}^{[i+1]}[\psi]$ at $(x_1,\dots,x_{n+2})$.
\end{proof}

\begin{remark}\label{rem:symmetric-reduction}
The hypotheses of Proposition~\ref{prop:symmetric-reduction} concern the datum alone, not the operator. When $\mathsf{T}$ is associative and symmetric, they hold on the image. Indeed, every $\psi\in\operatorname{Im}\mathsf{T}^{\,n}$ is a symmetric function of its $n+1$ arguments by Theorem~\ref{thm:symmetric-towers}, so the $n$ identities of Corollary~\ref{cor:higher-cocycle} carry no more information than the first. The necessary conditions obtained in this way for the inverse problem of $\mathsf{T}^{\,n}$ are thus the symmetry of the datum and one cocycle identity. For $\mathsf{T}=\mathsf{L}$ and $\mathsf{T}=\mathsf{D}_{\langle 2\rangle}$, both symmetric when $G$ is commutative, this applies to every power. For $\mathsf{T}=\mathsf{Q}$ and $\mathsf{T}=\mathsf{P}$, which are not symmetric when $2G\neq 0$ and $H$ is non-trivial, it applies to the data produced by the even functions. If $f$ is even, then $\mathsf{Q}^{\,n}[f]$ and $\mathsf{P}^{\,n}[f]$ are symmetric functions of their $n+1$ arguments. Indeed, each inner sum in (\ref{eq:Qn}) is then the sum of $f\bigl(\sum_{k\in S}\varepsilon_k x_k\bigr)$ over the sign vectors $\varepsilon\in\{\pm1\}^{S}$, now indexed by the whole of~$S$, with $\varepsilon$ and $-\varepsilon$ identified, and this expression is symmetric in the variables indexed by~$S$. For $\mathsf{P}^{\,n}$ this is the single term $S=\{1,\dots,n+1\}$.
\end{remark}

For $\mathsf{T}=\mathsf{L}$ and $n=2$ Corollary~\ref{cor:higher-cocycle} produces necessary conditions for the equation $\mathsf{F}[f]=\psi$ of the Fréchet operator. When $G$ is commutative the two identities it gives are, by Proposition~\ref{prop:symmetric-reduction}, equivalent on a symmetric datum to the single one $\mathsf{L}^{[1]}[\psi]=\mathsf{L}^{[2]}[\psi]$, and symmetry of the datum is itself necessary. For $\mathsf{T}=\mathsf{Q}$ Corollary~\ref{cor:higher-cocycle} gives necessary conditions for the higher-order inverse problem of Question~\ref{q:towers-Q}.

\subsection*{The factorisation principle}
Corollary~\ref{cor:higher-cocycle} and the direct implication of Proposition~\ref{prop:cocycle-complex} are instances of an elementary principle.

\begin{proposition}[factorisation principle]\label{prop:factor-principle}
Let $r,n,s\ge 1$, and let $\mathsf{U}\colon \mathcal{H}_r\to \mathcal{H}_n$ and $\mathsf{A},\mathsf{B}\colon \mathcal{H}_n\to \mathcal{H}_s$ be operators with $\mathsf{A}\,\mathsf{U} = \mathsf{B}\,\mathsf{U}$. Then every $\psi\in\operatorname{Im}\mathsf{U}$ satisfies $\mathsf{A}[\psi]=\mathsf{B}[\psi]$.
\end{proposition}

\begin{proof}
If $\psi=\mathsf{U}[\chi]$, then $\mathsf{A}[\psi]=\mathsf{A}\mathsf{U}[\chi]=\mathsf{B}\mathsf{U}[\chi]=\mathsf{B}[\psi]$.
\end{proof}

The lifts supply such pairs $(\mathsf{A},\mathsf{B})$. Associativity provides the pair $(\mathsf{T}^{[1]},\mathsf{T}^{[2]})$ for $\mathsf{U}=\mathsf{T}$, and Theorem~\ref{thm:towers} the pairs $(\mathsf{T}^{[i]},\mathsf{T}^{[j]})$ for $\mathsf{U}=\mathsf{T}^{\,n}$. The relation (\ref{eq:R}) provides the pair $(\mathsf{T}^{[i]},\mathsf{T}^{[i+1]})$ for $\mathsf{U}=\mathsf{T}^{[i]}$. Definition~\ref{def:coherent} covers the case in which the operator to be inverted and the operator that produces the conditions are different.

\begin{definition}\label{def:coherent}
An operator $\mathsf{U}\colon \mathcal{H}_1\to \mathcal{H}_n$ is said to be \emph{coherent with} $\mathsf{T}\colon \mathcal{H}_1\to \mathcal{H}_2$ if $\mathsf{T}^{[i]}\,\mathsf{U}$ does not depend on $i\in\{1,\dots,n\}$.
\end{definition}

\begin{corollary}\label{cor:coherent}
If $\mathsf{U}$ is coherent with $\mathsf{T}$, then every $\psi\in\operatorname{Im}\mathsf{U}$ satisfies
\[
\mathsf{T}^{[1]}[\psi] = \mathsf{T}^{[2]}[\psi] = \dots = \mathsf{T}^{[n]}[\psi].
\]
\end{corollary}

\begin{proof}
Apply Proposition~\ref{prop:factor-principle} with $\mathsf{A}=\mathsf{T}^{[i]}$ and $\mathsf{B}=\mathsf{T}^{[j]}$, which satisfy $\mathsf{A}\,\mathsf{U}=\mathsf{B}\,\mathsf{U}$ by Definition~\ref{def:coherent}.
\end{proof}

\begin{example}\label{ex:coherent}
An associative operator is coherent with itself, and $\mathsf{T}^{\,n}$ is coherent with $\mathsf{T}$, by Theorem~\ref{thm:towers}. Let now $G$ be a commutative monoid and $H$ an abelian group. The substitution rule of Lemma~\ref{lem:substitution} says that the mixed difference at the origin $\mathsf{D}_{\langle n\rangle}$ is coherent with $\mathsf{L}$, so that, by Corollary~\ref{cor:coherent}, every $\psi\in\operatorname{Im}\mathsf{D}_{\langle n\rangle}$ satisfies the cocycle identities
\[
\mathsf{L}^{[1]}[\psi] = \mathsf{L}^{[2]}[\psi] = \dots = \mathsf{L}^{[n]}[\psi],
\]
which are necessary conditions for the solvability of the inhomogeneous equation $\mathsf{D}_{\langle n\rangle}[f]=\psi$.
\end{example}

\subsection*{The image of a composition operator}
For composition operators the cocycle identity can also be sufficient. Over monoids and over bands it characterises the image, for an arbitrary target and with an explicit inverse, whereas over some semigroups it does not.

\begin{theorem}[monoids]\label{thm:monoid-image}
Let $G$ be a monoid with operation $\mu$ and identity element~$e$, and let $\psi\in \mathcal{H}_2$. Then $\psi\in\operatorname{Im}\mu^{*}$ if and only if, for all $x,y,z\in G$,
\begin{equation}\label{eq:sincov}
\psi\bigl(\mu(x,y),z\bigr) = \psi\bigl(x,\mu(y,z)\bigr).
\end{equation}
In that case $f := \psi(\,\cdot\,,e)$ is the unique solution of $\mu^{*}[f]=\psi$.
\end{theorem}

\begin{proof}
Throughout the proof we write $xy:=\mu(x,y)$. If $\psi=f\circ\mu$, both sides of (\ref{eq:sincov}) equal $f(xyz)$, and $f(x)=f(xe)=\psi(x,e)$, so the solution is unique. Conversely, the instance $(x,y,e)$ of (\ref{eq:sincov}) reads $\psi(xy,e)=\psi(x,ye)=\psi(x,y)$, so $\psi=\mu^{*}[\psi(\,\cdot\,,e)]$.
\end{proof}

\begin{proposition}[bands]\label{prop:bands}
Let $G$ be a band with operation $\mu$, and let $\psi\in \mathcal{H}_2$. Then $\psi\in\operatorname{Im}\mu^{*}$ if and only if $\psi$ satisfies (\ref{eq:sincov}). In that case $f(t):=\psi(t,t)$ is the unique solution of $\mu^{*}[f]=\psi$.
\end{proposition}

\begin{proof}
Throughout the proof we write $xy:=\mu(x,y)$. The operation $\mu$ is associative, hence so is $\mu^{*}$ by Proposition~\ref{prop:model}, and necessity is Proposition~\ref{prop:cocycle-complex}. A second solution $u$ satisfies $u(t)=u(tt)=\psi(t,t)$, so the solution is unique. Conversely, let $\psi$ satisfy (\ref{eq:sincov}). Its instances $(x,x,y)$ and $(x,xy,xy)$ read
\[
\begin{gathered}
\psi(x,y)=\psi(xx,y)=\psi(x,xy),\\
\psi(xy,xy)=\psi(x(xy),xy)=\psi(x,(xy)(xy))=\psi(x,xy),
\end{gathered}
\]
by $xx=x$, $x(xy)=(xx)y=xy$ and $(xy)(xy)=xy$. Hence $\psi(x,y)=\psi(xy,xy)=f(xy)$ for all $x,y\in G$.
\end{proof}

For a linearly ordered set under $\max$, which is a semilattice, Proposition~\ref{prop:bands} applies with $f(t)=\psi(t,t)$. Over an arbitrary semigroup, however, identity (\ref{eq:sincov}) need not characterise the image. Let $G$ be a \emph{null semigroup} with at least two elements, that is, a semigroup with an element $\theta$ such that $\mu(x,y)=\theta$ for all $x,y\in G$. Then (\ref{eq:sincov}) reads $\psi(\theta,z)=\psi(x,\theta)$. It requires $\psi$ to take one and the same value at all points of $\{\theta\}\times G$ and of $G\times\{\theta\}$, and nothing else, whereas $\operatorname{Im}\mu^{*}$ consists of the constant functions. If $H$ has at least two elements, a function that is constant on these two sets and takes a different value at some $(x,x)$ with $x\neq\theta$ satisfies (\ref{eq:sincov}) and is not in the image. The general criterion uses the following construction.

Let $G$ be a semigroup with operation $\mu$. Identity (\ref{eq:sincov}) says that $\psi$ is \emph{balanced}, that is, that it factors through the quotient of $G^2$ by the equivalence relation $\sim$ generated by $(\mu(x,y),z)\sim(x,\mu(y,z))$. We denote this quotient by $G\otimes_G G$, in the notation of the tensor product of acts \cite[§II.5]{KilpKnauerMikhalev}. The condition $\psi\in\operatorname{Im}\mu^{*}$ says that $\psi$ factors through~$\mu$. Since $\mu$ is associative, it is constant on the classes of~$\sim$ and induces a map $\bar\mu\colon G\otimes_G G\to G$, and the difference between the two conditions is measured by~$\bar\mu$. In elementary terms, $\bar\mu$ is injective when any two pairs with the same value under $\mu$ are joined by a chain of moves $(\mu(x,y),z)\leftrightarrow(x,\mu(y,z))$.

\begin{proposition}[balancedness]\label{prop:balanced}
Let $G$ be a semigroup with operation $\mu$. If $\bar\mu$ is injective, then a function $\psi\in \mathcal{H}_2$ lies in $\operatorname{Im}\mu^{*}$ if and only if it satisfies (\ref{eq:sincov}). If $\bar\mu$ is not injective and $H$ has at least two elements, then some $\psi\in \mathcal{H}_2$ satisfies (\ref{eq:sincov}) and does not lie in $\operatorname{Im}\mu^{*}$.
\end{proposition}

\begin{proof}
A function $\psi\in \mathcal{H}_2$ satisfies (\ref{eq:sincov}) precisely when it is constant on the classes of~$\sim$, that is, when there is a function $\bar\psi$ on $G\otimes_G G$ with $\psi(x,y)=\bar\psi(p)$ for $p$ the class of $(x,y)$. For the same $p$ one has $\mu(x,y)=\bar\mu(p)$, by the definition of~$\bar\mu$.

Every $\psi\in\operatorname{Im}\mu^{*}$ satisfies (\ref{eq:sincov}), by Proposition~\ref{prop:cocycle-complex}. Conversely, let $\bar\mu$ be injective and let $\psi$ satisfy (\ref{eq:sincov}), with $\bar\psi$ the function it determines. Since $\bar\mu$ is injective, $f\bigl(\bar\mu(p)\bigr):=\bar\psi(p)$ is a well-defined function on $\operatorname{Im}\bar\mu$. Extend it to $G$ arbitrarily. Then $\mu^{*}[f](x,y)=f\bigl(\bar\mu(p)\bigr)=\bar\psi(p)=\psi(x,y)$ for every pair $(x,y)$, with $p$ its class.

Let now $\bar\mu(p)=\bar\mu(q)$ for two classes $p\neq q$, and let $H$ contain $c_0\neq c_1$. Let $\bar\psi$ take the value $c_0$ at $p$ and the value $c_1$ everywhere else, and let $\psi$ be the function it defines on $G^2$, which satisfies (\ref{eq:sincov}). Every class is the class of a pair, so $\psi=\mu^{*}[f]$ would give $\bar\psi=f\circ\bar\mu$, whence $c_0=\bar\psi(p)=f\bigl(\bar\mu(p)\bigr)=f\bigl(\bar\mu(q)\bigr)=\bar\psi(q)=c_1$, a contradiction.
\end{proof}

Over a monoid the proof of Theorem~\ref{thm:monoid-image} gives the relation $(x,y)\sim(\mu(x,y),e)$, and over a band the proof of Proposition~\ref{prop:bands} gives $(x,y)\sim(\mu(x,y),\mu(x,y))$. In both cases $\bar\mu$ is bijective. A semigroup with $\bar\mu$ bijective is called \emph{firm} in the Morita theory of \cite{LaanMarkiReimaa2018}, and every semigroup with weak local units is firm \cite[Proposition~2.4]{LaanMarkiReimaa2018}. Monoids and bands have weak local units, so Theorem~\ref{thm:monoid-image} and Proposition~\ref{prop:bands} express their firmness for composition operators, with the inverse made explicit. Firmness is sufficient for (\ref{eq:sincov}) to characterise the image, but it is not necessary, since Proposition~\ref{prop:balanced} requires only injectivity. Consider the positive integers under addition. The instance $(x-1,1,y)$ of (\ref{eq:sincov}) identifies $(x,y)$ with $(x-1,y+1)$ whenever $x\ge 2$, so the classes of $\sim$ are the fibres of the sum and $\bar\mu$ is injective. It is not surjective, because $1$ is not a sum of two positive integers, so this semigroup is not firm. Yet (\ref{eq:sincov}) characterises $\operatorname{Im}\mu^{*}$ there, and the solutions are given by $f(t)=\psi(t-1,1)$ for $t\ge 2$, with $f(1)$ arbitrary. For a null semigroup with at least two elements, $\bar\mu$ is not injective.

We turn to the weighted composition operators $\omega\mu^{*}$ of Section~\ref{sec:lifts}, with weights in a commutative ring $R$ and with an $R$-module $H$ as target. For a weight $\omega\colon G^2\to R^{\times}$ with invertible values and for $\psi\in \mathcal{H}_2$ we put $\psi_{\omega}(x,y):=\omega(x,y)^{-1}\psi(x,y)$. The identity
\begin{equation}\label{eq:weighted-sincov}
\omega(x,y)\,\psi\bigl(\mu(x,y),z\bigr) = \omega(y,z)\,\psi\bigl(x,\mu(y,z)\bigr),
\end{equation}
required for all $x,y,z\in G$, is the weighted form of (\ref{eq:sincov}). When $\mu$ is associative and $\omega$ satisfies (\ref{eq:mult-cocycle}), it is the cocycle identity (\ref{eq:cocycle-cond}) of the operator $\omega\mu^{*}$. An invertible weight that satisfies (\ref{eq:mult-cocycle}) can be divided out, which reduces the weighted problem to the unweighted one.

\begin{lemma}[division by the weight]\label{lem:weights}
Let $G$ be a semigroup with operation $\mu$, let $R$ be a commutative ring, let $H$ be an $R$-module, and let $\omega\colon G^2\to R^{\times}$ satisfy (\ref{eq:mult-cocycle}). Let $\psi\in \mathcal{H}_2$ and $f\in \mathcal{H}_1$. Then
\begin{enumerate}
\item[\rm (i)] $(\omega\mu^{*})[f]=\psi$ if and only if $\mu^{*}[f]=\psi_{\omega}$;
\item[\rm (ii)] $\psi$ satisfies (\ref{eq:weighted-sincov}) if and only if $\psi_{\omega}$ satisfies (\ref{eq:sincov}).
\end{enumerate}
\end{lemma}

\begin{proof}
Since the values of $\omega$ are invertible, the equality $\psi(x,y)=\omega(x,y)f(\mu(x,y))$ holds if and only if $\psi_{\omega}(x,y)=f(\mu(x,y))$, which is~(i). For~(ii), substituting $\psi=\omega\psi_{\omega}$ in (\ref{eq:weighted-sincov}) gives
\[
\omega(x,y)\,\omega\bigl(\mu(x,y),z\bigr)\,\psi_{\omega}\bigl(\mu(x,y),z\bigr)
= \omega(y,z)\,\omega\bigl(x,\mu(y,z)\bigr)\,\psi_{\omega}\bigl(x,\mu(y,z)\bigr),
\]
and the two coefficients coincide by (\ref{eq:mult-cocycle}) and are invertible, so cancelling them leaves (\ref{eq:sincov}) for~$\psi_{\omega}$.
\end{proof}

\begin{corollary}[weighted monoids and bands]\label{cor:weighted}
Under the hypotheses of Lemma~\ref{lem:weights}, let $G$ be a monoid with identity element~$e$, or a band. Then the equation $(\omega\mu^{*})[f]=\psi$ is solvable if and only if $\psi$ satisfies (\ref{eq:weighted-sincov}). The solution is then unique, namely $f(t)=\omega(t,e)^{-1}\psi(t,e)$ for a monoid and $f(t)=\omega(t,t)^{-1}\psi(t,t)$ for a band.
\end{corollary}

\begin{proof}
By Lemma~\ref{lem:weights} the equation is equivalent to $\mu^{*}[f]=\psi_{\omega}$ and the condition to (\ref{eq:sincov}) for~$\psi_{\omega}$. Theorem~\ref{thm:monoid-image} and Proposition~\ref{prop:bands} then give solvability together with the unique solutions $f=\psi_{\omega}(\,\cdot\,,e)$ and $f(t)=\psi_{\omega}(t,t)$, which are the two formulas stated.
\end{proof}

\subsection*{Slicing and the stratification of towers}
For $\psi\in \mathcal{H}_{n+1}$, $1\le i\le n$ and a choice of the passive variables, the function of the active pair
\[
(s,t) \longmapsto \psi(x_1,\dots,x_{i-1},s,t,x_{i+2},\dots,x_{n+1})
\]
is called a \emph{slice} of $\psi$ in the slots $i$ and~$i+1$. The inverse problem for a lift reduces to the inverse problem for the operator itself, slice by slice.

\begin{proposition}[slicing]\label{prop:slicing}
Let $\mathsf{T}\colon \mathcal{H}_1\to \mathcal{H}_2$ be an operator, let $n\ge 1$, $1\le i\le n$ and $\psi\in \mathcal{H}_{n+1}$. Then $\psi\in\operatorname{Im}\mathsf{T}^{[i]}$ if and only if every slice of $\psi$ in the slots $i$ and $i+1$ lies in $\operatorname{Im}\mathsf{T}$.
\end{proposition}

\begin{proof}
By (\ref{eq:lift}), the equation $\mathsf{T}^{[i]}[\chi]=\psi$ splits along the passive tuples: for each tuple~$p$ it requires exactly $\mathsf{T}[\chi_p]=\psi_p$, where $\chi_p$ and $\psi_p$ are the corresponding partial maps. Distinct passive tuples involve disjoint families of values of~$\chi$, so solutions $u_p$ of the slice equations, chosen arbitrarily, assemble into $\chi(x_1,\dots,x_n) := u_p(x_i)$, and conversely.
\end{proof}

\begin{remark}[stratification]\label{rem:stratify}
Let $\mathsf{T}$ be associative, $n\ge 2$ and $\psi\in \mathcal{H}_{n+1}$. Then $\psi\in\operatorname{Im}\mathsf{T}^{\,n}$ if and only if there exists $\chi\in\operatorname{Im}\mathsf{T}^{\,n-1}$ with $\mathsf{T}^{[n]}[\chi]=\psi$: take $\chi=\mathsf{T}^{\,n-1}[f]$ in (\ref{eq:power-rec}) with $i=n$, and conversely $\mathsf{T}^{[n]}\mathsf{T}^{\,n-1}[f]=\mathsf{T}^{\,n}[f]$. In particular, by Proposition~\ref{prop:slicing}, every slice of $\psi$ in the last two variables lies in $\operatorname{Im}\mathsf{T}$. The inverse problem of order~$n$ is thus reduced to the inverse problem for $\mathsf{T}$ itself, applied slice by slice, together with a selection problem one level down: among the solutions of the slice equations, one must choose a family that assembles into an element of $\operatorname{Im}\mathsf{T}^{\,n-1}$. This description of $\operatorname{Im}\mathsf{T}^{\,n}$ by means of $\operatorname{Im}\mathsf{T}^{\,n-1}$ is called the \emph{stratification} of the towers of~$\mathsf{T}$.
\end{remark}

Beyond the cocycle identity the answer depends on the operator. For $\mathsf{L}$ over a commutative semigroup symmetry must be added, and the two conditions together are sufficient over a cancellative abelian monoid with a divisible target \cite{DavisonEbanks1995}. For $\mathsf{Q}$ symmetry is not available when $2G\neq 0$ and $H$ is non-trivial, and the quasi-symmetry of Remark~\ref{rem:quasi-symmetry} takes its place among the necessary conditions. For a composition operator over a general semigroup the question is the injectivity of $\bar\mu$ (Proposition~\ref{prop:balanced}). Question~\ref{q:towers-Q} asks for the corresponding conditions for the higher powers $\mathsf{Q}^{\,n}$. Where the inverse is explicit, as in Theorem~\ref{thm:monoid-image} and Proposition~\ref{prop:bands}, each value of the solution is a value of the datum, so a bound on the datum is a bound on the solution. This is the connection with stability mentioned in Section~\ref{sec:intro}.

\section{New versus known and open questions}\label{sec:new}

For composition operators the calculus of lifts is classical, and this paper claims no novelty there. For operators $\mu^{*}$ induced by an associative binary operation, the lifts are the duals of the maps that fold two adjacent coordinates, and the relation (\ref{eq:R}) is a simplicial coface relation. In operadic terms it follows from the associativity of $\mu$ and the sequential composition axiom for partial compositions, while the exchange law of Remark~\ref{rem:no-exchange} is the parallel composition axiom \cite[§5.3]{LodayVallette}. The multiplace operations $G^n\to G$, with the partial compositions, form the endomorphism operad of~$G$, and the algebras of multiplace functions with Menger's superposition are treated in \cite{DudekTrokhimenko}. For composition operators Theorem~\ref{thm:towers} reduces to the generalised associative law \cite[Ch.~I, §1, no.~3]{Bourbaki}, and the full set of relations is given by the universal property of the simplicial category of \cite[VII.5, Proposition~1]{MacLaneCWM}, whose objects are all the finite ordinals, the empty one included.

The image question for a composition operator is known as well, in a different language. Balancedness is the notion on which the tensor product of acts is built \cite[§II.5]{KilpKnauerMikhalev}, and the bijectivity of $\bar\mu$ is the firmness of the Morita theory of \cite{LaanMarkiReimaa2018}, which the discussion after Proposition~\ref{prop:balanced} relates to Theorem~\ref{thm:monoid-image} and Proposition~\ref{prop:bands}. Proposition~\ref{prop:balanced} follows directly from these notions. We have not found it, nor the identification of the cocycle identity of a composition operator with balancedness, in the literature on functional equations.

The higher-order Cauchy differences and their kernels are studied in \cite{NgZhao2013,GuoMaLi,Ng2015,LiNg2016,Stetkaer2017}, and the compatibility of two associative products, in the bilinear setting, is studied in \cite{OdesskiiSokolov}. For the addition of an abelian group, condition (\ref{eq:mult-cocycle}) on the weight is the classical associativity condition of a twisted group ring \cite[Ch.~1, §2]{Passman1977}, as Remark~\ref{rem:twisted} records. Proposition~\ref{prop:L}, the associativity of the Cauchy difference, is the classical cocycle identity of that operator, read through Definition~\ref{def:associative}. Proposition~\ref{prop:Q}, the associativity of $\mathsf{Q}$, is in the same way the identity established in \cite{PopaRusJvN}, where the inverse problem for $\mathsf{Q}$ is studied. The quasi-symmetry of Remark~\ref{rem:quasi-symmetry} is an elementary expansion, recorded here for comparison and not claimed as new. The solution theory of the quadratic equation recalled in Section~\ref{sec:lifts}, and the theory of the cocycle equation of the Cauchy difference recalled in Section~\ref{sec:image}, are classical as well.

What we believe to be new is the following. The notion of an associative operator is new, together with the lifts in every slot and the tower theorem (Theorem~\ref{thm:towers}). The two-slot identity that defines the notion is anticipated in \cite{PopaRusJvN}, as Section~\ref{sec:lifts} records. Example~\ref{ex:no-exchange} shows that the notion differs from the cosimplicial one: an associative operator satisfies the relation (\ref{eq:R}) in each slot, but its lifts need not satisfy the exchange law. In Section~\ref{sec:sums}, the classification of Theorem~\ref{thm:global} and the resulting characterisations of $\mathsf{L}$, $\mathsf{P}$, $\mathsf{Q}$ and $\mathsf{D}_{\langle 2\rangle}$ (Corollary~\ref{cor:detect}) are new, and so is the operator form of compatibility (Proposition~\ref{prop:sum}), with its equivalent form in Remark~\ref{rem:compatible-combinations} and the two compatible pairs of Corollary~\ref{cor:pencils}. So are the associativity of the d'Alembert operator, the closed form of the powers of $\mathsf{Q}$, the description of the towers of $\Delta$ as the partial mixed differences, the associativity of the integer combinations of $\mathsf{L}$ and $\mathsf{J}_0$, and the identification of the mixed differences at the origin, from the second one on, as the powers of $\mathsf{L}+\mathsf{J}_0$ (Propositions \ref{prop:P}, \ref{prop:Qn}, \ref{prop:delta-towers}, \ref{prop:D2} and~\ref{prop:origin-tower}). The results on symmetry of Section~\ref{sec:symmetry} are new, except Proposition~\ref{prop:symmetric-model}, which is elementary. In Section~\ref{sec:image}, Proposition~\ref{prop:cocycle-complex} is Definition~\ref{def:associative} read on the image. It gives the identities for $\mathsf{L}$, for $\mathsf{Q}$ and for the composition operators a common form, which \cite{PopaRusJvN} observes for the first two. New in that section are the extension of the identity to the powers (Corollary~\ref{cor:higher-cocycle}), its reduction to a single identity on symmetric data (Proposition~\ref{prop:symmetric-reduction}), the notion of coherence (Definition~\ref{def:coherent}), the slicing of lifts (Proposition~\ref{prop:slicing}) and the stratification of towers (Remark~\ref{rem:stratify}). The factorisation principle of Proposition~\ref{prop:factor-principle} is only the observation that two operators which agree after precomposition with a third agree on its image. Lemma~\ref{lem:weights} is an elementary reduction, and Corollary~\ref{cor:weighted} follows from it and from the unweighted results. The tower spectrum of Definition~\ref{def:tower-spectrum} extends the associative spectrum from composition operators to arbitrary operators. It is new in that generality, and so are the values of Proposition~\ref{prop:tower-spectrum}.

The exchange law, the classification over the targets excluded by the hypothesis of Theorem~\ref{thm:global}, and the conditions on the image that lie beyond the cocycle identity are taken up in a sequel in preparation.

We close with three questions.

\begin{question}\label{q:towers-Q}
The powers $\mathsf{Q}^{\,n}$ are given in closed form by (\ref{eq:Qn}). Does the inverse problem for $\mathsf{Q}^{\,n}$, in the spirit of \cite{PopaRusJvN} for $n=1$, admit a uniform solution? Corollary~\ref{cor:higher-cocycle} supplies necessary conditions at every order. What replaces the quasi-symmetry of Remark~\ref{rem:quasi-symmetry} for $n\ge 2$?
\end{question}

\begin{question}\label{q:weighted-bands}
Lemma~\ref{lem:weights} divides out a weight that takes invertible values and satisfies (\ref{eq:mult-cocycle}), and so reduces the weighted image problem to the unweighted one, where over an arbitrary semigroup the question is the injectivity of $\bar\mu$ (Proposition~\ref{prop:balanced}). Both hypotheses on the weight are used. If $\omega$ takes non-invertible values, the weight cannot be divided out. If $\omega$ fails (\ref{eq:mult-cocycle}), the operator $\omega\mu^{*}$ need not be associative, and then Proposition~\ref{prop:cocycle-complex} gives no condition. What characterises $\operatorname{Im}\omega\mu^{*}$ in these two cases?
\end{question}

The last question concerns operators that are not associative, through the tower spectrum of Definition~\ref{def:tower-spectrum}.

\begin{question}\label{q:tower-spectrum}
Which sequences occur as tower spectra?
\end{question}

\section*{Statements and Declarations}

\subsection*{Funding}

No funding was received to assist with the preparation of this manuscript.

\subsection*{Competing Interests}

The authors have no competing interests to declare that are relevant to the
content of this article.

\subsection*{Use of AI Tools}

The authors used the large language model Claude (Anthropic) in the preparation of this manuscript: to check the results and their proofs, to search and check the literature, and to revise the text. The authors have verified all the statements, proofs and references, and take full responsibility for the content of the article.

\bibliographystyle{spmpsci}
\bibliography{references}

\end{document}